\documentclass[10pt, final, journal, letterpaper, twoside, twocolumn]{IEEEtran}
\usepackage{amsmath,amssymb}
\usepackage{enumerate}
\usepackage{graphicx}
\usepackage{float}
\usepackage{stackengine}
\usepackage[caption=false,font=footnotesize]{subfig}
\usepackage{diagbox}
\usepackage[noadjust]{cite}
\usepackage{multirow}
\usepackage{booktabs}
\usepackage{adjustbox}
\usepackage{color}
\renewcommand{\Re}{\mathbb{R}}
\newcommand{\Rc}{\overline{\Re}}
\newcommand{\Lin}{\mathcal{L}}

\newcommand{\mcR}{\mathcal{R}}

\newcommand*{\rom}[1]{\expandafter\@slowromancap\romannumeral #1@}

\def\D{\mathbf{D}}
\def\W{\mathbf{W}}
\def\A{\mathbf{A}}
\def\B{\mathbf{B}}
\def\D{\mathbf{D}}

\def\I{\mathbf{I}}
\def\K{\mathbf{K}}
\def\W{\mathbf{W}}
\def\U{\mathbf{U}}
\def\V{\mathbf{V}}
\def\P{\mathbf{P}}
\def\W{\mathbf{W}}

\def\ZE{\mathbf{0}}

\def\c{\boldsymbol{c}}
\def\e{\boldsymbol{e}}
\def\h{\boldsymbol{h}}
\def\i{\boldsymbol{i}}

\def\n{\boldsymbol{n}}

\def\r{\boldsymbol{r}}
\def\s{\boldsymbol{s}}
\def\t{\boldsymbol{t}}
\def\u{\boldsymbol{u}}
\def\v{\boldsymbol{v}}

\def\x{\boldsymbol{x}}
\def\y{\boldsymbol{y}}
\def\z{\boldsymbol{z}}

\def\prox{\mathrm{Prox}}

\def\fix{\mathrm{fix}}

\DeclareMathOperator*{\argmin}{\arg\!\min}

\newcommand{\inner}[2]{\left\langle#1,#2\right\rangle}
\newcommand{\innerE}[2]{\inner{#1}{#2}_2}
\newcommand{\norm}[1]{\left\lVert#1\right\rVert}
\newcommand{\normE}[1]{\norm{#1}_2}

\graphicspath{{./figures/}}

\newtheorem{theorem}{Theorem}[section]
\newtheorem{lemma}[theorem]{Lemma}
\newtheorem{proposition}[theorem]{Proposition}

\newtheorem{definition}[theorem]{Definition}

\begin{document}

% ---- Title -----------------------------------------
\title{Fixed-Point and Objective Convergence of Plug-and-Play Algorithms}
% ----------------------------------------------------

% ---- Authors ---------------------------------------
\author{Pravin Nair,~\IEEEmembership{Student~Member,~IEEE}, Ruturaj G. Gavaskar,~and~Kunal~N.~Chaudhury,~\IEEEmembership{Senior~Member,~IEEE}
% ----------------------------------------------------

% ---- Affiliations ---------------------------------
\thanks{The authors are with the Department of Electrical Engineering, Indian Institute of Science, Bengaluru, India. The work of K.~N.~Chaudhury was supported by MATRICS grant MTR/2018/000221 from the Department of Science and Technology, Government of India, and Grant ISTC/EEE/KNC/440 from the ISRO-IISc Space Technology Cell.}}

\markboth{Submitted to IEEE Transactions on Computational Imaging}{}
\maketitle

\begin{abstract}
A standard model for image reconstruction involves the minimization of a data-fidelity term along with a regularizer, where the optimization is performed using proximal algorithms such as ISTA and ADMM. In plug-and-play (PnP) regularization, the proximal operator (associated with the regularizer) in ISTA and ADMM  is replaced by a powerful image denoiser. Although PnP regularization works surprisingly well in practice, its theoretical convergence---whether convergence of the PnP iterates is guaranteed and if they minimize some objective function---is not completely understood even for simple linear denoisers such as nonlocal means. In particular, while there are works where either iterate or objective convergence is established separately, a simultaneous guarantee on iterate and objective convergence is not available for any denoiser to our knowledge. In this paper, we establish both forms of convergence for a special class of linear denoisers. Notably, unlike existing works where the focus is on symmetric denoisers, our analysis covers non-symmetric denoisers such as nonlocal means and almost any convex data-fidelity. The novelty in this regard is that we make use of the convergence theory of averaged operators and we work with a special inner product (and norm) derived from the linear denoiser; the latter requires us to appropriately define the gradient and proximal operators associated with the data-fidelity term. We validate our convergence results using image reconstruction experiments. 
\end{abstract}


\section{Introduction}

In imaging modalities such as tomography and MRI, we are required to reconstruct a high-resolution  image $\x_0 \in \Re^n$ from incomplete noisy measurements $\y \in \Re^m$  \cite{ribes2008linear,scherzer2009variational}. Whereas, in applications such as single-image superresolution and deblurring, a subsampled or blurred image $\y$ is available and we need to infer the ground-truth from $\y$ and some knowledge of the degradation process \cite{dong2011image}. More generally, the abstract problem of inverting a given measurement model comes up in several  computational imaging applications \cite{scherzer2009variational,engl1996regularization}. The standard optimization framework for addressing such problems involves a data-fidelity term $f: \Re^n \to \Re$ and a regularizer $g: \Re^n \to \Re \cup \{\infty\}$; the former is derived from the measurement model, while the latter is typically derived from Bayesian or sparsity-promoting priors \cite{scherzer2009variational}. The reconstruction is given by the solution of the  optimization problem 
\begin{equation}
\label{mainopt}
\underset{\x \in \Re^n}{\min} \  f(\x) + g(\x).
\end{equation}
For example, in deblurring, superresolution, and compressive imaging \cite{dong2011image,jagatap2019algorithmic}, the measurement model is linear and $f$ is given by $f(\x) = \lambda \lVert \A\x - \y \rVert^2$, where $\A \in \mathbb{Re}^{m \times n}$ is the measurement matrix and $\lambda > 0$ is used to balance the terms in \eqref{mainopt} (in this paper, we absorb $\lambda$ in the data-fidelity term for convenience). 
On the other hand, the data-fidelity term is non-quadratic but convex  in problems such as single photon imaging \cite{chan2017plug}, Poisson denoising \cite{rond2016poisson}, and despeckling \cite{bioucas2010multiplicative}. 

%\subsection{Regularization using denoising} 

The regularizer $g$ in (1) penalizes signals that look significantly different from natural images; this in effect forces the reconstruction to resemble the ground-truth. From computational considerations, $g$ is often taken to be convex. In fact, there has been significant amount of research on the design of convex regularizers \cite{ribes2008linear,scherzer2009variational,engl1996regularization}. More recently, it was shown in several works that off-the-shelf denoisers can be used for regularization within an iterative framework. For example, in Plug-and-Play (PnP) regularization \cite{sreehari2016plug}, this is done by fixing a proximal algorithm (henceforth referred to as the \textit{base algorithm}) for solving \eqref{mainopt} and replacing the proximal operator of $g$ with a powerful denoiser. In particular, this has been done for proximal algorithms such as iterative shrinkage thresholding algorithm (ISTA) \cite{Sun2019_PnP_SGD,Ryu2019_PnP_trained_conv}, and alternating direction method of multipliers (ADMM) \cite{sreehari2016plug,Ryu2019_PnP_trained_conv}.
Remarkably, PnP regularization (or simply `PnP') has been shown to produce promising results across a wide range of applications including superresolution, MRI, fusion, and tomography 
\cite{sreehari2016plug,zhang2017learning,Dong2018_DNN_prior,Sun2019_PnP_SGD,Ryu2019_PnP_trained_conv,
rick2017one,zhang2019deep,kamilov2017plug,ono2017primal, Tirer2019_iter_denoising,Teodoro2019_PnP_fusion,song2020new,ahmad2020plug}. We note that denoisers have also been used for regularization using different schemes \cite{bigdeli2017deep,romano2017little,reehorst2018regularization,mataev2019deepred,sun2019block}.
%Henceforth, we will specifically use the term ``PnP regularization'' (or just PnP) when a denoiser is used within a base algorithm for regularization.

\subsection{Motivation} 

Although PnP works well in practice, there is apriori no reason why the PnP iterates should converge in the first place.
Moreover, it is not clear whether they minimize some objective of the form in \eqref{mainopt}. The former, commonly referred to as {\em fixed-point convergence}, was investigated in \cite{chan2017plug,Gavaskar2019_fixed_point,Ryu2019_PnP_trained_conv}.
On the other hand, the question of optimality and {\em objective convergence} was addressed for a class of symmetric linear denoisers in \cite{sreehari2016plug,Teodoro2019_PnP_fusion}.
The convergence guarantee in these existing works hold for specific classes of denoisers.
For example, the denoiser is assumed to be bounded in \cite{chan2017plug,Tirer2019_iter_denoising}, symmetric in \cite{Teodoro2019_PnP_fusion}, averaged in \cite{Sun2019_PnP_SGD}, demicontractive in \cite{cohen2020regularization}, and the residue corresponding to the denoiser is assumed to be non-expansive in \cite{Ryu2019_PnP_trained_conv} (some of these terms will be defined later in the paper).
In practice, though, it is  difficult to verify if a denoiser is bounded or averaged. 
Along with the denoiser, the data-fidelity is often constrained as well. 
For example, the convergence guarantee in \cite{Ryu2019_PnP_trained_conv} requires the data-fidelity to be strongly convex,  which is not met for imaging problems like superresolution, compressed sensing, despeckling, etc.

The focus of this work is on linear denoisers and, in particular, kernel filters such as the Yaroslavsky filter \cite{yaroslavsky1985digital}, Lee filter \cite{lee1983digital}, bilateral filter \cite{tomasi1998bilateral}, nonlocal means (NLM) \cite{buades2005non}, LARK \cite{takeda2007kernel} etc.
In particular, PnP using NLM-type denoisers is known to produce promising image reconstructions \cite{venkatakrishnan2013plug,heide2014flexisp,sreehari2016plug,sreehari2017multi,unni2018linearized,Chan2019_PnP_graph_SP, nair2019hyperspectral,gavaskar2020plug,unni2020_pnp_registration}.
Since linear operators are easier to analyze, establishing PnP convergence for linear denoisers is a natural first step towards understanding the behavior of nonlinear denoisers.
Even for linear denoisers, existing theoretical guarantees are limited to symmetric denoisers \cite{sreehari2016plug,Teodoro2019_PnP_fusion}; notably, this excludes neighborhood filters such as NLM that are naturally non-symmetric \cite{Singer2009_diffusion_nonlocal,Milanfar2013_filtering_tour}. 

\subsection{Prior work}
\label{prior}

The problem of PnP convergence has been studied in several works.
It was shown in \cite{sreehari2016plug} that a class of symmetric linear denoisers can be expressed as the proximal operator of a convex function, i.e., one can associate a convex regularizer $g$ for every denoiser in this class, and the PnP iterates in this case amount to minimizing $f+g$. The closed-form expression of $g$ was later derived in \cite{Teodoro2019_PnP_fusion}.
However, it is generally difficult to certify whether a denoiser can be expressed as a proximal map; this is particularly true for complex denoisers such as BM3D \cite{dabov2007image}, TNRD \cite{chen2016trainable} and DnCNN \cite{zhang2017beyond}. 
The next best in this case is to guarantee convergence of the PnP iterates. This has been done for different combinations of inverse problem, base algorithm and denoiser. For example, iterate convergence was established for linear inverse problems with quadratic data-fidelity in \cite{Dong2018_DNN_prior,Tirer2019_iter_denoising,gavaskar2020plug}. In addition, the denoiser is assumed to satisfy a descent condition in \cite{Dong2018_DNN_prior}, a boundedness condition in \cite{Tirer2019_iter_denoising}, and a linearity condition in \cite{gavaskar2020plug}. On the other hand, the analysis in \cite{chan2017plug,Sun2019_PnP_SGD} applies to arbitrary convex data-fidelity but  the denoiser is  assumed to satisfy a boundedness condition in \cite{chan2017plug} (similar to \cite{Tirer2019_iter_denoising}), an averagedness property in \cite{Sun2019_PnP_SGD} and demicontractivity in \cite{cohen2020regularization}. 

The convergence of PnP algorithms using deep denoisers has been studied in recent papers.
It was shown in \cite{Ryu2019_PnP_trained_conv} that PnP convergence is guaranteed for ISTA and ADMM for a specially trained CNN denoiser, provided the data-fidelity is strongly convex.
Apart from \cite{Ryu2019_PnP_trained_conv}, PnP convergence has been established for CNN denoisers \cite{Meinhardt2017_learning_prox_op,buzzard2018plug}, generative denoisers \cite{jagatap2019algorithmic}, and GAN-based projectors \cite{raj2019gan}.
Moreover, it was shown in \cite{xu2020provable} that the DnCNN denoiser can be approximately expressed as the proximal operator of a nonconvex function.

\subsection{Contribution}

In this work, we establish iterate and objective convergence for a class of non-symmetric linear denoisers. Importantly, we do not assume the data-fidelity to be strongly convex. On the technical front, our key contributions are as follows: 

\begin{enumerate}

\item We prove that convergence of the PnP iterates can be guaranteed for ISTA and ADMM provided $f$ is convex and the denoiser is averaged (see Definition \ref{def:averagedop}); $f$ must be differentiable for ISTA but is allowed to be non-smooth for ADMM. Our analysis is based on the theory of proximal and averaged operators \cite{bauschke2017convex}. In particular, we prove the existence of non-symmetric linear denoisers that are averaged with respect to a non-standard norm.

\item In Theorems \ref{thm:linearISTA} and \ref{thm:linearADMM}, we simultaneously establish  iterate and objective convergence for a special class of non-symmetric denoisers. In particular, this subsumes existing results on symmetric denoisers \cite{sreehari2016plug,Teodoro2019_PnP_fusion}. Notably, unlike \cite{Ryu2019_PnP_trained_conv,raj2019gan,fletcher2018plug}, our results hold for arbitrary convex data-fidelity. Our analysis highlights the need to work with a non-standard inner product derived from the denoiser. In effect, this requires us to appropriately define the gradient and proximal operators in ISTA and ADMM. 

\item In Theorem \ref{thm:NLMkernelfilt}, we prove that the NLM denoiser is the proximal map of a quadratic convex function, where the norm in the proximal map is induced by a non-standard inner product.

\end{enumerate}

We validate our theoretical findings and demonstrate the effectiveness of the modified PnP algorithms using superresolution and despeckling experiments. 

\subsection{Organization}
In Section \ref{background}, we collect some basic notations, definitions and algorithmic specifications. Our results on PnP convergence are stated and discussed  in Section \ref{mainresults}; detailed derivations of these results are deferred to Section \ref{sec:proofs}. We validate our findings using superresolution and despeckling experiments in Section \ref{Experiments} and conclude  with a discussion in Section \ref{Conclusion}. 

\section{Preliminaries}
\label{background}

\subsection{Notations}

We denote the standard inner product on $\Re^n$ by $\innerE{\cdot}{\cdot}$, i.e.,
$\innerE{\x}{\y} = \x^\top \!\y$ for $\x,\y \in \Re^n$.
 The norm induced by $\innerE{\cdot}{\cdot}$ is denoted by $\normE{\cdot}$.
We denote the identity operator on $\Re^n$ by $I$, whereas the identity matrix is denoted by $\I$.
The range space of a matrix $\A$ is denoted by $\mcR(\A)$.  

We will work with non-standard inner products and norms on $\Re^n$. In particular, we will require the notion of gradient, non-expansivity, proximal operator, etc. in a real Hilbert space $(\Re^n, \inner{\cdot}{\cdot})$, where $\inner{\cdot}{\cdot}$ is an abstract inner product on $\Re^n$. We let $\norm{\cdot}$ denote the norm induced by $\inner{\cdot}{\cdot}$, i.e., $\norm{\x}=\inner{\x}{\x}^{1/2}$ for $\x \in \Re^n$. We will define the PnP iterations using the gradient and proximal operators corresponding to $\inner{\cdot}{\cdot}$.

\subsection{Basic Definitions}
\label{Prelims}

We begin by defining non-expansive and averaged operators on $\Re^n$.  
\begin{definition}
\label{def:averagedop}
%Let $\inner{\cdot}{\cdot}$ be an inner product on $\Re^n$, and $\norm{\cdot}$ be its induced norm.
An operator $T$ on $(\Re^n, \inner{\cdot}{\cdot})$ is said to be non-expansive  if $\norm{T (\x) - T(\y)} \leqslant \norm{\x - \y}$ for all $\x,\y \in \Re^n$.
An operator $A$ on $(\Re^n, \inner{\cdot}{\cdot})$ is said to be $\theta$-averaged, where $\theta \in (0,1),$ if we can write $A = (1-\theta) I + \theta T$, where $T$ is non-expansive. That is, $A$ is $\theta$-averaged if $(1 - 1/\theta)I + (1/\theta)A$ is non-expansive.
\end{definition}
%Equivalently, $T$ is $\theta$-averaged if $(1 - 1/\theta)I + (1/\theta)T$ is non-expansive.

We next define the gradient and the proximal operator in a Hilbert space $(\Re^n, \inner{\cdot}{\cdot})$.
\begin{definition}
\label{def:gradient}
A function $f : \Re^n \to \Re$ is said to be differentiable at $\x \in \Re^n$ if there exists a (unique) linear map $L : \Re^n \to \Re$ (the derivative of $f$ at $\x$) such that
\begin{equation}
\label{diff}
f(\x+\h) = f(\x) + L(\h) + o(\norm{\h}) \qquad (\h \in \Re^n),
\end{equation}
where $\norm{\cdot}$ is any arbitrary norm on $\Re^n$. If \eqref{diff} holds at every $\x \in \Re^n$, then we say that $f$ is differentiable and we use $L(\x)$ to denote the derivative at $\x$. For a fixed inner product $\inner{\cdot}{\cdot}$ on $\Re^n$, the corresponding gradient of $f$ at $\x$ is the unique vector $\nabla \! f(\x) \in \Re^n$ such that
\begin{equation}
\label{grad}
L(\x)(\h) = \inner{\nabla \! f(\x)}{\h} \qquad (\h \in \Re^n).
\end{equation}

\end{definition}

Note that the usual definition of $\nabla \! f$ as the vector of partial derivatives of $f$ is consistent with Definition~\ref{def:gradient} when the inner product is $\innerE{\cdot}{\cdot}$ \cite{rudin1964principles}.

The proximal operator \cite{bauschke2017convex} is at the heart of algorithms such as ISTA and ADMM \cite{Parikh2014proximal}. Let $\Rc$ denote the extended real line $\Re \cup \{\infty\}$. An extended-real-valued  $f : \Re^n \to \Rc$ is said to be proper if there exists $\x \in \Re^n$ such that $f(\x) < \infty$.  Moreover, $f$ is said to be closed if its epigraph,
\begin{equation*}
\mathrm{epi}(f) = \big\{(\x,t) \in \Re^{n+1}: t \geqslant f(\x) \big\},
\end{equation*}
is closed in $\Re^{n+1}$. 

\begin{definition}
\label{def:prox_op}
Let $f : \Re^n \to \Rc$ be a closed, proper and convex function. 
The proximal operator  $\prox_f$ on $(\Re^n, \inner{\cdot}{\cdot})$ is defined as
\begin{equation}
\label{prox}
\prox_f(\x) = \argmin_{\y \in \Re^n}\  \frac{1}{2} \norm{\y - \x}^2 + f(\y),
\end{equation}
where $\norm{\cdot}$ is the norm induced by $\inner{\cdot}{\cdot}$.
\end{definition}

\subsection{Algorithms}
The standard ISTA and ADMM algorithms for solving \eqref{mainopt} can be generalized to $(\Re^n, \inner{\cdot}{\cdot})$ (where the inner product $\inner{\cdot}{\cdot}$ is arbitrary), using the  gradient and proximal operators in \eqref{grad} and \eqref{prox}.
The ISTA  iterations in $(\Re^n, \inner{\cdot}{\cdot})$ are given by
\begin{align*}
&{\z}_{k+1} = {\x}_{k} -  \rho^{-1} \nabla \! f ({\x}_{k}). \\
&{\x}_{k+1} = \prox_{\rho^{-1} g} (\z_{k+1}). 
\end{align*}
where $\rho > 0$ and $\x_0 \in \Re^n$ is an arbitrary initialization. On the other hand, the ADMM iterations are given by
\begin{align*}
\x_{k+1} &= \prox_{{\rho}^{-1} f} \left(\y_k - \z_k \right), \\
\y_{k+1} &= \prox_{{\rho}^{-1} g} (\x_{k+1} + \z_k), \\
\z_{k+1} &= \z_k + \x_{k+1} - \y_{k+1}, 
\end{align*}
where $\rho > 0$ and $\y_0,\z_0 \in \Re^n$ are initializations. In the PnP variant of ISTA (ADMM), referred to as PnP-ISTA (PnP-ADMM), the proximal operator $\prox_{\rho^{-1} g}$ is replaced by an image denoiser $D$.  

\section{Convergence Analysis}
\label{mainresults}

In this section, we establish the following results for PnP-ISTA and PnP-ADMM:
\begin{itemize}
\item If the denoiser is averaged, then the iterates of PnP-ISTA and PnP-ADMM exhibit fixed-point convergence.
\item The  averaged property  (with respect to a special inner product) is satisfied by a broad class of linear denoisers.
\item For this class of linear denoisers, there exists a convex regularizer $g$ such that the limit points of PnP-ISTA and PnP-ADMM are minimizers of $f + g$. 
\end{itemize}

We will just state and discuss these results and connect them to existing results; their detailed derivations can be found in Section \ref{sec:proofs}. As is well-known, the proximal operator (see Definition~\ref{def:prox_op}) of a closed, proper, and convex function is $(1/2)$-averaged in $(\Re^n, \innerE{\cdot}{\cdot})$ \cite{Parikh2014proximal}. 
As a result, a denoiser that can be expressed as the proximal operator of a convex regularizer is averaged.
For this reason, the symmetric linear denoisers in \cite{sreehari2016plug,Teodoro2019_PnP_fusion} qualify as averaged operators on $(\Re^n, \innerE{\cdot}{\cdot})$. But what about a generic linear denoiser $D(\x) = \W \x$, where $\W$ has the basic properties of nonnegativity and row-stochasticity, but is possibly non-symmetric?
%The following result highlights the importance of working with an inner product different from $ \innerE{\cdot}{\cdot}$ in this case.
The following result highlights that such denoisers do not qualify as averaged operators.

\begin{proposition}
\label{prop:averaged_l2}
Let $D$ be a linear operator on $\Re^n$. In particular, let $D(\x) = \W \x$, where $\W \in \Re^{n \times n}$.
\begin{enumerate}[(a)]
\item Suppose $\W$ is symmetric and its eigenvalues are in $[0,1]$. Then $D$ is $\theta$-averaged on  $(\Re^n, \innerE{\cdot}{\cdot})$ for all $\theta \in [1/2,1)$.
\item Let $\W$ be nonnegative and row-stochastic (all rows sum to one), but not doubly-stochastic (some columns do not sum to one). Then $D$ cannot be $\theta$-averaged on $(\Re^n, \innerE{\cdot}{\cdot})$ for any $\theta \in (0,1)$.
\end{enumerate}
\end{proposition}

In particular, Proposition~\ref{prop:averaged_l2} implies that kernel filters, such as those mentioned earlier \cite{yaroslavsky1985digital,lee1983digital,tomasi1998bilateral,buades2005non,takeda2007kernel} are not averaged on $(\Re^n, \innerE{\cdot}{\cdot})$.
We remark that the weight matrix $\W$ is derived from the input image $\x$ in most of these filters.
As a result, these filters are not strictly linear in terms of the input-output relation.
However, $\W$ can be computed from a surrogate image \cite{Milanfar2013_filtering_tour} or fixed after a few PnP iterations \cite{heide2014flexisp,sreehari2016plug,sreehari2017multi,unni2018linearized,Chan2019_PnP_graph_SP,
nair2019hyperspectral,gavaskar2020plug}. The filter can be treated as a linear operator $D(\x) = \W \x$ in this case \cite{Singer2009_diffusion_nonlocal}. In particular, while $\W$ is nonnegative and row-stochastic, it is naturally non-symmetric. Hence, $D$ cannot possibly be $\theta$-averaged on $(\Re^n, \innerE{\cdot}{\cdot})$.

The above negative result leads us to the natural question: Is $\W$ averaged with respect to some non-standard inner product?
We show in Proposition~\ref{prop:neighborfilter} that this is indeed the case. In particular, we will fix an appropriate inner product $\inner{\cdot}{\cdot}$ on $\Re^n$ and  use the corresponding gradient and proximal operators within PnP-ISTA and PnP-ADMM.  More precisely, we consider the following PnP-ISTA iterations:
\begin{equation}
\label{eq:PnP_ISTA}
{\x}_{k+1} = D ({\x}_{k} - \rho^{-1} \nabla \! f ({\x}_{k})),
\end{equation}
and the following PnP-ADMM iterations:
\begin{subequations}
\label{eq:PnP_ADMM}
\begin{align}
\x_{k+1} &= \prox_{{\rho}^{-1} f} \left(\y_k - \z_k \right), \label{eq:PnP-ADMM_x}\\
\y_{k+1} &= D (\x_{k+1} + \z_k), \label{eq:PnP-ADMM_y}\\
\z_{k+1} &= \z_k + \x_{k+1} - \y_{k+1}. \label{eq:PnP-ADMM_z}
\end{align}
\end{subequations}
With the averaged property in place for kernel filters, we can establish convergence of  PnP-ISTA and PnP-ADMM using such filters. Notably, we do not need to symmetrize $\W$ that entails additional cost \cite{sreehari2016plug}.

We now analyze the fixed-point convergence of PnP-ISTA and PnP-ADMM for averaged denoisers $D : \Re^n \to \Re^n$.
This is based on the fixed-point theory of averaged operators \cite{bauschke2017convex}.
\begin{definition}
\label{def:fixed_point}
We say that $\x^\ast \in \Re^n$ is a fixed point of $T: \Re^n \to \Re^n$ if $T(\x^\ast) = \x^\ast$.
The set of fixed points of $T$ is denoted by $\mathrm{fix}(T)$.
\end{definition}

We state a lemma \cite{bauschke2017convex} that is required to establish convergence of PnP-ISTA and PnP-ADMM.
\begin{lemma}
\label{lem:key_lemma}
Let $T :\Re^n \to \Re^n$ be $\theta$-averaged on $(\Re^n, \inner{\cdot}{\cdot}),$ where $\theta \in (0,1)$. 
Assume that $\fix(T)$ is not empty and let  $\x_0 \in \Re^n$. Then the sequence $(\x_k)_{k \geqslant 0}$ generated as $\x_{k+1} = T(\x_k)$ converges to some $\x^\ast \in \fix(T)$.
\end{lemma}

%A well-known result for averaged operators is that they exhibit fixed-point convergence 
To apply this result to PnP-ISTA, we need the notion of a smooth function.
\begin{definition}
\label{def:smooth_fn}
Let $f : \Re^n \to \Re$ be differentiable. It is said to be $\beta$-smooth on $(\Re^n, \inner{\cdot}{\cdot})$ if there exists $\beta >0$ such that $\norm{\nabla \! f(\x) - \nabla \! f(\y)} \leqslant \beta \norm{\x - \y}$ for all $\x, \y \in \Re^n$, where $\norm{\cdot}$ is the norm induced by $\inner{\cdot}{\cdot}$.
\end{definition}

We are now ready to state our main results on the fixed-point convergence of PnP-ISTA and PnP-ADMM.
Henceforth, we assume that the data-fidelity term $f : \Re^n \to \Re$ is real-valued (rather than extended real-valued), which is the case with most imaging applications.
\begin{theorem}%[Fixed-point convergence of PnP-ISTA]
\label{thm:PnP_ISTA_fixed}
%Consider the PnP-ISTA algorithm on $(\Re^n, \inner{\cdot}{\cdot})$ given by \eqref{eq:PnP_ISTA}.
Let $T_\mathrm{ISTA} = D \circ (I - \rho^{-1} \nabla \! f)$.
Suppose that
\begin{itemize}
\item $f$ is convex and $\beta$-smooth,
\item $D$ is $\theta$-averaged for some $\theta \in (0,1)$, and 
\item $\fix(T_\mathrm{ISTA}) \neq \varnothing$. 
\end{itemize}
Then for any $\x_0 \in \Re^n$ and $\rho > \beta/2$, the sequence $(\x_k)_{k \geqslant 0}$ generated by \eqref{eq:PnP_ISTA} converges to some $\x^\ast \in \fix(T_\mathrm{ISTA})$.
\end{theorem}
We remark that fixed-point convergence of PnP-ISTA for a larger class of denoisers, including averaged denoisers, was recently established in \cite{cohen2020regularization}, although under slightly stricter assumptions. 

\begin{theorem}%[Fixed-point convergence of PnP-ADMM]
\label{thm:PnP_ADMM_fixed}
%Consider the PnP-ADMM algorithm on $(\Re^n, \inner{\cdot}{\cdot})$ given by \eqref{eq:PnP_ADMM}.
Define
\begin{equation}
\label{eq:T_ADMM}
T_{\mathrm{ADMM}} = \frac{1}{2} I + \frac{1}{2} \left(2 \prox_{\rho^{-1} f} - I \right) \circ \left(2 D - I \right).
\end{equation}
Suppose that
\begin{itemize}
\item $f$ is convex,
\item $D$ is $\theta$ averaged for some $\theta \in (0,1/2]$, and 
\item $\fix(T_{\mathrm{ADMM}}) \neq \varnothing$. 
\end{itemize}
Then, for arbitrary $\rho>0$ and $\y_0, \z_0 \in \Re^n$, the sequence $(\x_{k}, \y_{k}, \z_{k})_{k \geqslant 1}$ generated by \eqref{eq:PnP_ADMM} is convergent and the limit point is determined by some $\u^{\ast} \in \fix(T_{\mathrm{ADMM}})$.
\end{theorem}

Note that the results in Theorems \ref{thm:PnP_ISTA_fixed} and \ref{thm:PnP_ADMM_fixed} hold for any choice of the inner product $\inner{\cdot}{\cdot}$.
Moreover, $D$ is not assumed to be linear.
However, it is assumed that $\fix(T_{\mathrm{ISTA}})$ and $\fix(T_{\mathrm{ADMM}})$ are non-empty. 
There are a couple of issues in this regard.
First, verifying whether a given denoiser is averaged is not an easy task; this is especially true for non-linear denoisers such as BM3D \cite{dabov2007image} and DnCNN \cite{zhang2017beyond}.
Second, even if the denoiser $D$ is averaged, it is unclear whether $\fix(T_{\mathrm{ISTA}})$ and $\fix(T_{\mathrm{ADMM}})$ are non-empty.
In many cases, the latter condition is not verifiable and is simply assumed to hold without proof \cite{Sun2019_PnP_SGD}. 
In the subsequent discussion, we show how to deal with these issues for a special class of linear denoisers.

\begin{definition}
Let $\Lin$ denote the class of linear denoisers $D(\x) = \W \x$ on $\Re^n$ such that $\W$ is diagonalizable and its eigenvalues are in $[0,1]$.
\end{definition}

In the following theorem, we collect some relevant properties of $\Lin$, particularly that every denoiser  in $\Lin$ is averaged.
\begin{theorem}
\label{thm:averaged_combined}
Let $D$ be in class $\Lin$ and $\W$ be the associated matrix, i.e., $D(\x)=\W\x$. Let $\V \in \Re^{n \times n}$ be an eigen matrix of $\W$, i.e., the columns of $\V$ are linearly independent eigenvectors of $\W$. Define the inner product
\begin{equation}
\label{inner}
\inner{\x}{\y} = \innerE{\V^{-1} \x}{\V^{-1} \y}.
\end{equation}
Then we have the following properties:
\begin{enumerate}[(a)]
\item $D$ is the proximal operator on $(\Re^n,\inner{\cdot}{\cdot})$ of some closed, proper and convex function $g_{D}: \Re^n \to \overline{\Re}$. 
\item $D$ is $\theta$-averaged  on $(\Re^n,\inner{\cdot}{\cdot})$ for every $\theta \in [1/2,1)$.
\item The restriction of $g_{D}$ to $\mcR(\W)$ is real-valued (and hence continuous).  
\item The iterates of PnP-ISTA in \eqref{eq:PnP_ISTA} (resp. PnP-ADMM in \eqref{eq:PnP_ADMM}) are identical to that of ISTA (resp. ADMM) applied to the optimization problem 
\begin{equation}
\underset{\x \in \Re^n}{\min} \  f(\x) + \rho g_{D}(\x).
\label{eq:mainoptlinear}
\end{equation}
\end{enumerate}
\end{theorem}

Since \eqref{inner} depends on $D$ via an eigen basis of $\W$, we will refer to this as an inner product induced by $D$. 
We can interpret \eqref{inner} as the standard inner product applied along with a change of basis, i.e., we perform our computations  with respect to an eigen basis of $\W$ instead of the standard basis. In particular, if $\W$ is symmetric, then \eqref{inner} is in fact $\innerE{\cdot}{\cdot}$ since the eigen matrix $\V$ can be taken to be orthogonal in this case. 

If the denoiser $D$ is in class $\Lin$ and $\inner{\cdot}{\cdot}$ is an inner product induced by $D$, then the iterates in PnP-ISTA and PnP-ADMM  are guaranteed to converge to a minimizer of \eqref{eq:mainoptlinear}.
\begin{theorem}
Let $D \in \Lin$ and $\inner{\cdot}{\cdot}$ be an inner product induced by $D$.
Consider the space $(\Re^n,\inner{\cdot}{\cdot})$.
Assume that
\begin{itemize}
\item $f$ is convex and $\beta$-smooth.
\item $\x^\ast$ is a minimizer of \eqref{eq:mainoptlinear} and $p^\ast= f(\x^\ast) + \rho g_{D}(\x^\ast)$.
\end{itemize}
Then, for any $\x_0 \in \Re^n$ and $\rho > \beta/2$,
\begin{enumerate}[a)]
\item the PnP-ISTA iterates ${(\x_k)}_{k \geqslant 0}$ generated by \eqref{eq:PnP_ISTA} converge to a minimizer of \eqref{eq:mainoptlinear}, and
\item $\mathrm{lim}_{k \to \infty} \ f(\x_k) + \rho g_{D}(\x_k) = p^\ast$.  
\end{enumerate}
\label{thm:linearISTA}
\end{theorem}

\begin{theorem} Let $D \in \Lin$ and $\inner{\cdot}{\cdot}$ be an inner product induced by $D$. Consider the space $(\Re^n,\inner{\cdot}{\cdot})$.
Assume that $f$ is convex, $\y^\ast$ is a minimizer of \eqref{eq:mainoptlinear} and $p^\ast= f(\y^\ast) + \rho g_{D}(\y^\ast)$. Then, for any $\rho>0$ and initialization $\y_0, \z_0 \in \Re^n$,
\begin{enumerate}[a)]
\item the PnP-ADMM iterates ${(\y_{k})}_{k \geqslant 0}$ generated by \eqref{eq:PnP-ADMM_y} converge to a minimizer of \eqref{eq:mainoptlinear}, and
\item $\mathrm{lim}_{k \to \infty} f(\y_k) + \rho g_{D}(\y_k) = p^\ast$.
\end{enumerate}
\label{thm:linearADMM}
\end{theorem}

We note that Theorem \ref{thm:linearADMM} subsumes the objective convergence result in  \cite{sreehari2016plug}, where $\W$ is assumed to be symmetric. 
The above theorems show that convergence of PnP algorithms can be extended to a larger class of linear denoisers $\Lin$ including non-symmetric denoisers, provided we work with a special inner product.  

The practical utility of the convergence results is that many kernel filters belong to the class $\Lin$.
Though this is well-known \cite{Singer2009_diffusion_nonlocal},  we explain why this is so for completeness.
Let $\x \in \Re^N$ be the  vectorized input image, where $N$ is the number of pixels.
The elements of $\x$ are $\{\x_{\s}: \s \in \Omega \subset \mathbb{Z}^2\}$, where $\Omega$ is the support of the image, $\s \in \Omega$ is a pixel location, and $\x_{\s}$ is the corresponding intensity value.
Let $\boldsymbol{\zeta}_{\s} \in \Re^M$ represent some feature at pixel $\s$. The output of a generic kernel filter is given by 
\begin{equation}
\label{neighbour}
{(D(\x))}_{\s} = \frac{\sum_{\t \in \Omega} \phi(\boldsymbol{\zeta}_{\s}, \boldsymbol{\zeta}_{\t}) \x_{\t}}{\sum_{\t \in \Omega} \phi(\boldsymbol{\zeta}_{\s}, \boldsymbol{\zeta}_{\t})}, 
\end{equation}
where the kernel function $\phi: \Re^M \times \Re^M \to \Re$ is nonnegative and symmetric \cite{Morel2014}.
To obtain the matrix representation of \eqref{neighbour}, let $\mathfrak{R} =\{\boldsymbol{\xi}_1,\boldsymbol{\xi}_2,....,\boldsymbol{\xi}_N\}$ be some ordering of the elements in $ \big\{ \boldsymbol{\zeta}_{\s} \in \Re^{M} : \s \in \Omega \big\}$. That is, for every $\ell \in [1,N]$, $\boldsymbol{\xi}_{\ell} = \boldsymbol{\zeta}_{\s}$ for some $\s \in \Omega$. For $i, j \in \{1, 2, \ldots N\}$, define the kernel matrix $\K \in \mathbb{R}^{N \times N}$  to be
\begin{equation}
\label{kernelmatrix}
\K_{i,j} = \phi(\boldsymbol{\xi}_i, \boldsymbol{\xi}_j),
\end{equation}
and the diagonal normalization matrix $\D \in \Re^{N \times N}$ as
\begin{equation}
\label{normmatrix}
\D_{i,i} = \sum_{j=1}^N \K_{i,j}. 
\end{equation}
We can then write \eqref{neighbour} as
\begin{equation}
\label{neigbrfilt}
D(\x) = \D^{-1} \K \x. 
\end{equation}  

\begin{definition}
\label{def:pdkernel}
A kernel function $\phi: \Re^M \times \Re^M \to \Re$ is said to be positive definite if for any $N \geqslant 1$, $c_1, \ldots c_N \in \Re$, and $\x_1, \ldots \x_N \in \Re^M$, 
\begin{equation*}
\sum_{i=1}^N \sum_{j=1}^N c_i c_j \phi(\x_i, \x_j) \geqslant 0. 
\end{equation*}
\end{definition}
\begin{proposition}
\label{prop:neighborfilter}
If the kernel $\phi$ in \eqref{neighbour} is positive definite and $D$ is given by \eqref{neigbrfilt}, then $D$ belongs to $\Lin$.
Furthermore, $D$ is $\theta$-averaged on $(\Re^n,\inner{\cdot}{\cdot})$  for every $\theta \in [1/2,1)$, where $\inner{\cdot}{\cdot}$ is given by
\begin{equation*}
\inner{\x}{\y} = \x^\top\! \D \y,% \qquad \norm{\v} = \sqrt{{\inner{\v}{\v}}}, 
\end{equation*}
with $\D$ as in \eqref{normmatrix}. 
\end{proposition}

For more information on kernel filters, we refer the reader to \cite{milanfar2013symmetrizing,Morel2014}. 
In our experiments, we use the nonlocal means (NLM) denoiser, which is a special instance of \eqref{neighbour}. 
In NLM, the feature vector is given by $\boldsymbol{\zeta}_{\s} = (\s, \P_{\s})$, where $\s \in \Re^{2}$ is the pixel coordinates and $\P_{\s}$ is a (vectorized) image patch around pixel $\s$ of a guide image (which can be different from the input image). 
%Let $\P_{\s} \in \Re^{N_p^2}$ be patch of size $N_p \times N_p$ around pixel at locations $\s \in \Omega$. The feature vector corresponding to a pixel $\s \in \Omega$ is $\y_s = (\s, \P_{\s})$, where $\y_s \in \Re^M$ with $M = N_p^2 + 2$. 
In particular, we consider the following NLM kernel
\begin{equation}
\label{nlmkernel}
\phi (\boldsymbol{\zeta}_{\s}, \boldsymbol{\zeta}_{\t} ) = \Lambda(\s - \t)  \kappa(\P_{\s} - \P_{\t}), 
\end{equation}
where $\kappa$ is Gaussian and the hat function $\Lambda: \Re^2 \to \Re$ is given by
\begin{equation*}
\Lambda (\s) = \prod_{i=1}^2 \left(1 - \Big | \frac{s_i}{N_s} \Big| \right)_+, 
\end{equation*}
where $(t)_+ = \max(0,t)$ and $N_s$ is the search radius. 
% and range kernel $\kappa: \Re^{{N_p}^2} \to \Re$ is a gaussian kernel $\mathcal{N}(\textbf{0},\sigma^2 \I)$. 
The kernel in \eqref{nlmkernel} is positive definite \cite{sreehari2016plug}. 
Importantly, if we define $\K$ and $\D$ as in \eqref{kernelmatrix} and \eqref{normmatrix} and the NLM denoiser using \eqref{neigbrfilt} with a fixed guide image (used to compute the kernel in \eqref{nlmkernel}), then the denoiser belongs to class $\Lin$. 
Furthermore, the denoiser is the proximal map of a quadratic convex regularizer. 
\begin{theorem}
\label{thm:NLMkernelfilt}
The kernel matrix $\K$ in NLM is positive definite. Furthermore the corresponding denoiser \eqref{neigbrfilt} is the proximal map on $(\Re^n,\inner{\cdot}{\cdot})$ of the convex regularizer
\begin{equation*}
g_{D}(\x) = \frac{1}{2} \x^\top \D (\K^{-1} \D - \I) \x, 
\end{equation*}
where $\inner{\cdot}{\cdot}$ is as specified in Proposition \ref{prop:neighborfilter}. 
\end{theorem}
%We note that by working with PnP algorithms in the specific Hilbert space, we avoid the need for symmetrizing the filter as in \cite{sreehari2016plug}. 
%Instead we can directly work with neighbourhood filters in its given form (provided the kernel is positive definite) and still prove convergence of PnP algorithms. 

We conclude this section with a discussion of the scope of our analysis.
Our proofs are intricately tied to the form of the updates in ISTA and ADMM; it is not clear whether they can be adapted to other PnP algorithms.
In particular, the methods we have used to prove convergence do not apply to accelerated variants such as PnP-FISTA \cite{Sun2019_PnP_SGD}.
Another limitation is that our results are restricted to linear kernel filters and cannot to applied to nonlinear denoisers. In particular, for a kernel filter to be treated as a linear denoiser, we are required to fix $\W$ after a finite number of PnP iterations; that is, we cannot adapt $\W$ using the reconstruction beyond a point.
%While our analysis is applicable to any denoiser belonging to the class $\mathcal{L}$, PnP convergence for symmetric denoisers belonging to $\mathcal{L}$ has been separately studied before \cite{sreehari2016plug,Teodoro2019_PnP_fusion}.
%As for non-symmetric denoisers belonging to $\mathcal{L}$, kernel filters are the only practical examples.
%Finally, we again recall an important difference between `standard' PnP algorithms and the PnP algorithms studied in the current work: to calculate $\nabla f$ in PnP-ISTA and $\prox_{ \rho^{-1} f}$ in PnP-ADMM, the appropriate (possibly non-standard) inner product must be used.
%Thus, unlike standard PnP methods, this involves more than just replacing $\prox_{{\rho}^{-1} g}$ by a denoiser.
%$\nabla f$ can be computed from the standard gradient of $f$, whereas numerical methods may be needed to compute $\prox_{ \rho^{-1} f}$.
%In rare cases, it may be possible to directly compute the latter; examples in this regard are shown in the next section.

\section{Numerical results}
\label{Experiments}

In this section, we validate the convergence results for PnP-ISTA and PnP-ADMM using a couple of image reconstruction experiments---superresolution and despeckling.
The data-fidelity $f$ is quadratic for the former and non-quadratic for the latter.
Importantly, $f$ is convex but not strongly convex.
%Thus, we are able to numerically validate our results under very general conditions.
%We choose these applications as both have non-strongly convex data-fidelity terms $f$, with $f$ being quadratic for superresolution and non-quadratic for despeckling.
We use NLM as the denoiser (with the kernel in \eqref{nlmkernel}) and we work in the space defined by the inner product $\inner{\cdot}{\cdot}$ in Proposition \ref{prop:neighborfilter}.
The matrix $\D$ is computed from the image obtained after five PnP iterations, and the weight matrix $\W$ is kept fixed thereafter.
Thus, the linear denoiser $D(\x) = \W \x$ belongs to $\mathcal{L}$.
The purpose of the experiments is solely to demonstrate that iterate and objective convergence are indeed achieved in practical imaging problems, as predicted by Theorems \ref{thm:linearISTA} and \ref{thm:linearADMM}.
%Hence the denoiser used from iteration number $6$ in the algorithm belongs to $\Lin$ and convergence results follow from Theorem \ref{thm:linearISTA} and \ref{thm:linearADMM}. 
In particular, we do not claim that our reconstructions are superior to existing methods, including PnP with other denoisers.
Nevertheless, since PnP algorithms in $(\Re^n,\inner{\cdot}{\cdot})$ (for non-standard inner products) have not been used till date, we compare the reconstruction quality with PnP in $(\Re^n,\innerE{\cdot}{\cdot})$ (i.e. standard PnP).
This is done to confirm that a similar reconstruction quality is obtained regardless of the inner product used.
We stress that in $(\Re^n,\innerE{\cdot}{\cdot})$, convergence guarantees for PnP-ISTA and PnP-ADMM with NLM denoiser are not available in full generality.
Thus, working with the appropriate (denoiser-induced) inner product offers the advantage of a better understanding of the PnP mechanism (via an objective function), in addition to establishing convergence. 

\subsection{Superresolution} 

The observation model for superresolution is given by
\begin{equation}
\label{sr}
\y = \mathbf{S}\B\tilde{\x} + \boldsymbol{\eta},
\end{equation}
where $\tilde{\x} \in \mathbb{R}^n$ is the unknown high-resolution image, $\y \in \mathbb{R}^m$ is the observed low-resolution image ($m<n$), $\B \in \mathbb{R}^{n \times n}$ is a circulant matrix corresponding to a blur kernel $b$, and $\mathbf{S} \in \mathbb{R}^{m \times n}$ is a binary sampling matrix that decimates $\B\x$ by a factor $K$ \cite{chan2017plug}.
Note that the vectorized forms of the images are used in \eqref{sr}. 
For white Gaussian noise $\boldsymbol{\eta}$ (standard deviation $\sigma$), the data-fidelity term corresponding to the maximum likelihood estimate of $\tilde{\x}$ is given by
\begin{equation}
\label{loss}
f(\x) = \frac{1}{2} \normE{\y - \mathbf{S}\B\x}^2.
\end{equation}
This is not strongly convex since $\mathbf{S}$ has a non-trivial null space. 
\begin{table}[t!]
\centering
\caption{PSNR/SSIM values for superresolution of Set12 dataset \cite{zhang2017beyond} using PnP-ISTA.The  $\sigma$  values are on a scale of $0$ to $255$.}
\label{tab:superres}
\begin{tabular}{llcc}
\toprule
\multicolumn{2}{c}{\diagbox{Settings}{Method}} & PnP in $(\Re^n,\innerE{\cdot}{\cdot})$ & PnP in $(\Re^n,\inner{\cdot}{\cdot})$ \vspace{0.2mm}\\
\toprule
\multirow{2}{*}{$K=2$} & $\sigma=5$ & $28.15/0.830$ & $28.20/0.836$ \\
& $\sigma=10$ & $26.80/0.775$ & $26.94/0.785$ \\
\hline
\multirow{2}{*}{$K=4$} & $\sigma=5$ & $24.31/0.737$ & $24.29/0.726$ \\
& $\sigma=10$ & $23.38/0.689$ & $23.43/0.692$ \\
\bottomrule
\end{tabular}
\end{table}
We use PnP-ISTA to estimate the ground-truth high-resolution image. The gradient of \eqref{loss} in  $(\Re^n, \inner{\cdot}{\cdot})$ where $\inner{\cdot}{\cdot}$ is as specified in Proposition \ref{prop:neighborfilter}, is given by
\begin{equation*}
\nabla \! f(\x) = \D^{-1} \B^\top {\mathbf{S}}^\top\! (\mathbf{S}\B\x - \y).
\end{equation*}
Note that $\mathbf{S}^\top \in \mathbb{R}^{n \times m}$ is an upsampling matrix and $\B^\top \in \mathbb{R}^{n \times n}$ is again a circulant matrix whose blur kernel is obtained by flipping $b$ about the origin \cite{chan2017plug}.
In particular, $b$ is a symmetric Gaussian blur for our experiments and $\B^\top\! = \B$ in this case.
%\footnote{In practice, $\B^\top\!$ and $\B$ are never populated and stored. They are applied using FFT-based convolution.}. 
For all experiments, we use a symmetric Gaussian blur of size $9 \times 9$ and standard deviation $1$.

\begin{figure}
\centering
\scriptsize
\stackunder{\hphantom{=}}{\stackunder{\includegraphics[width=0.242\linewidth]{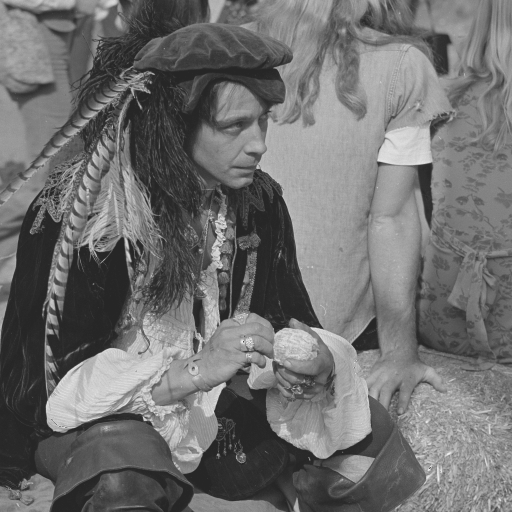}}{Ground truth.}}
\stackunder{\hphantom{=}}{\stackunder{\includegraphics[width=0.242\linewidth]{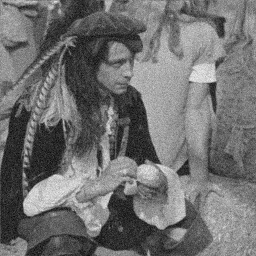}}{Observed.}}
\stackunder{PSNR = $27.72$ dB.}{\stackunder{\includegraphics[width=0.242\linewidth]{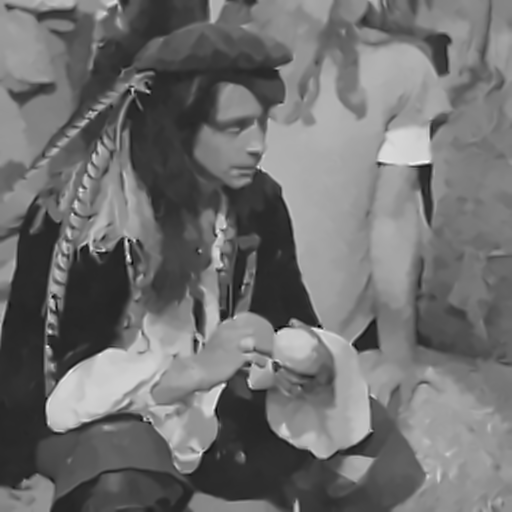}}{PnP, $(\Re^n,\innerE{\cdot}{\cdot})$.}}
\stackunder{PSNR = $28.01$ dB.}{\stackunder{\includegraphics[width=0.242\linewidth]{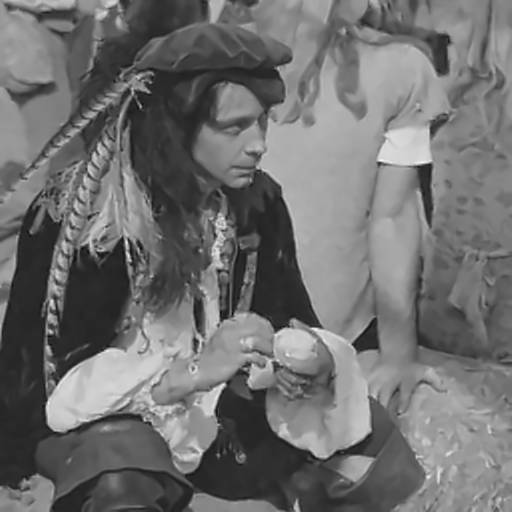}} {PnP, $(\Re^n,\inner{\cdot}{\cdot})$.}}
\caption{Image superresolution for $K=2$ and $\sigma = 10$ using PnP-ISTA. The observed image has size $256 \times 256$ (scaled for display purposes), whereas the other images have size $512 \times 512$. The third and fourth images are the reconstructions from PnP-ISTA in $(\Re^n,\innerE{\cdot}{\cdot})$ and $(\Re^n,\inner{\cdot}{\cdot})$. 
The respective SSIM values are $0.728$ and $0.747$, while the PSNR values are reported above the images.}
\label{fig:superres}
\end{figure}

\begin{figure}
\centering
\subfloat{\includegraphics[width=0.475\linewidth]{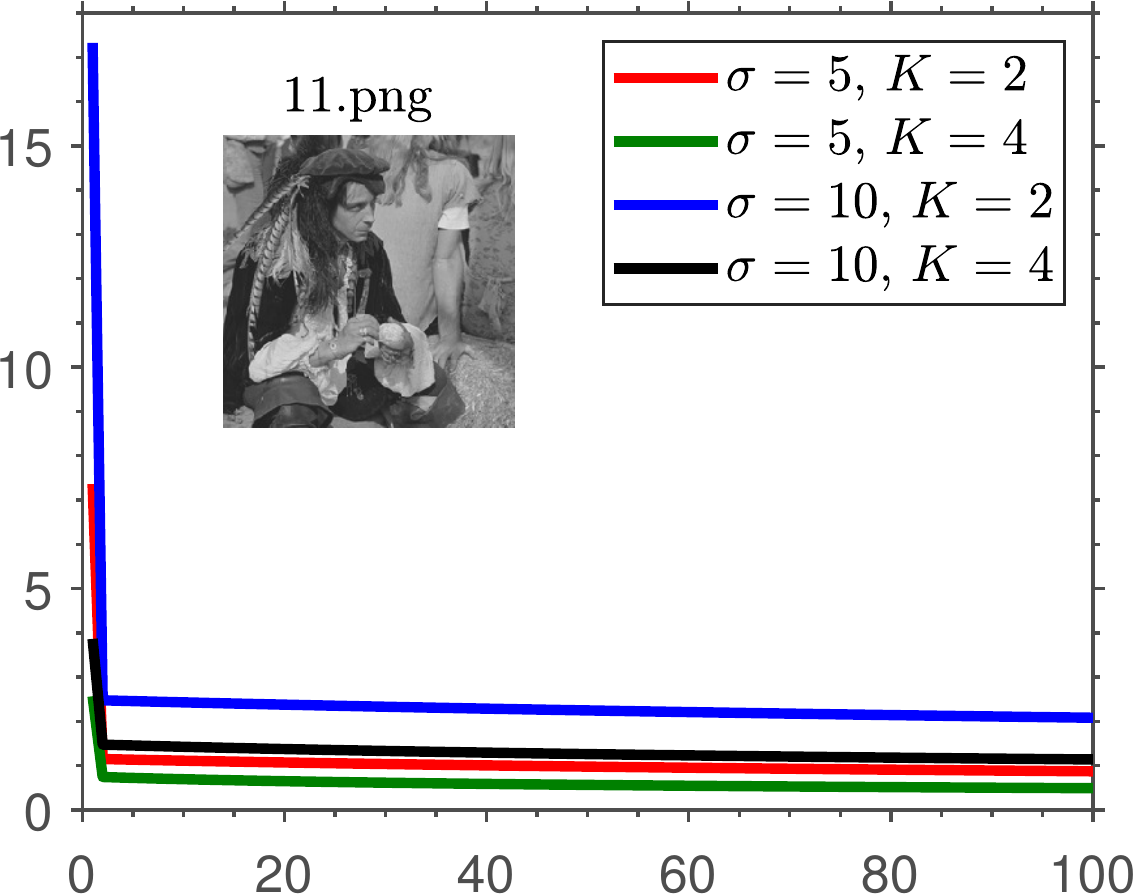}}  \hspace{1mm}
\subfloat{\includegraphics[width=0.475\linewidth]{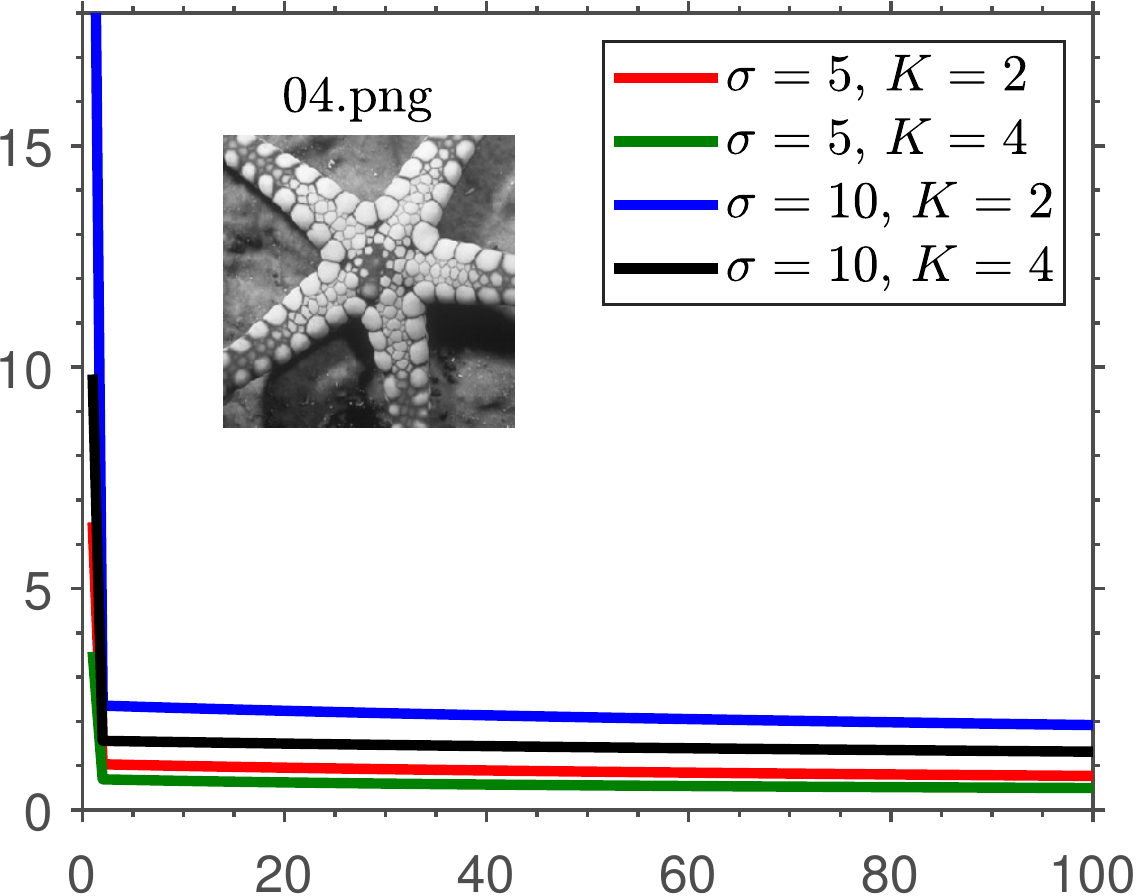}}  

\subfloat{\includegraphics[width=0.475\linewidth]{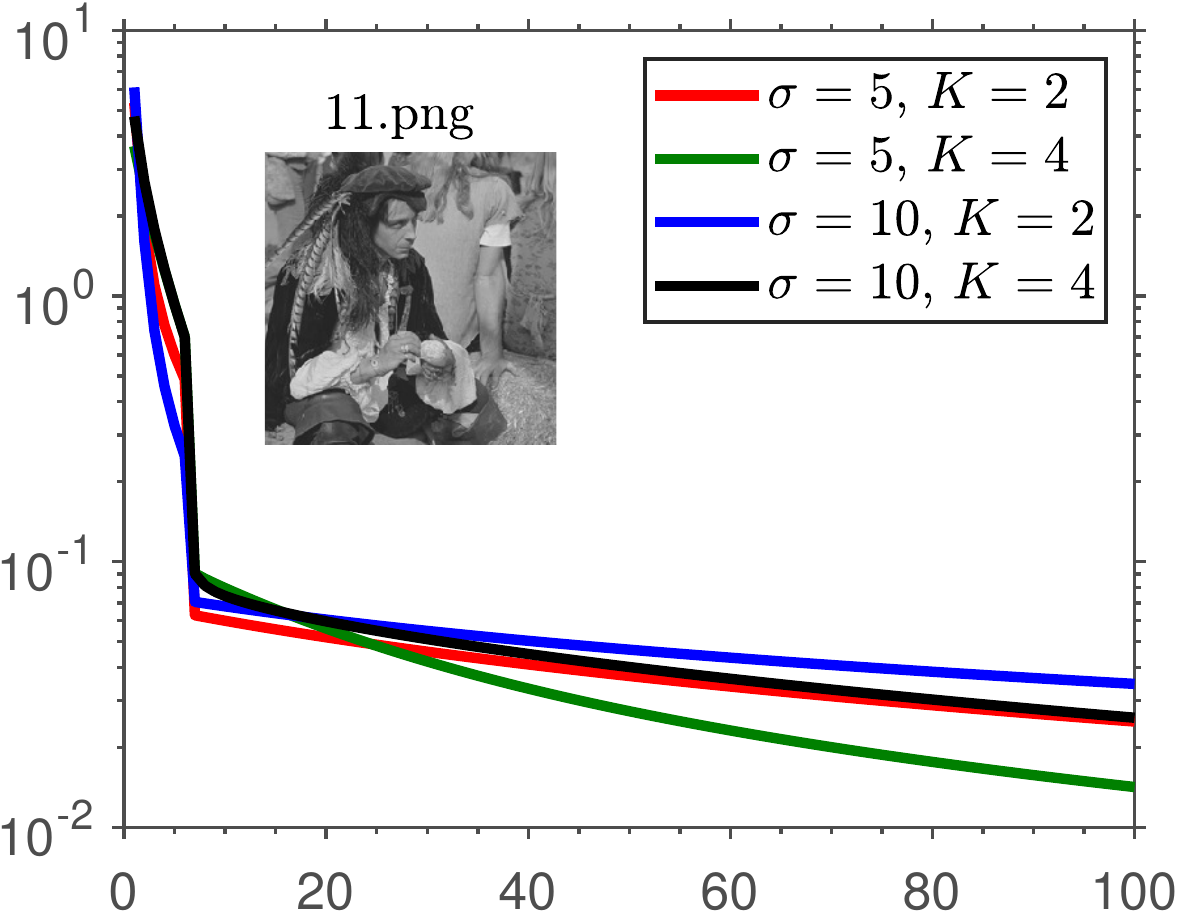}} \hspace{1mm}
\subfloat{\includegraphics[width=0.475\linewidth]{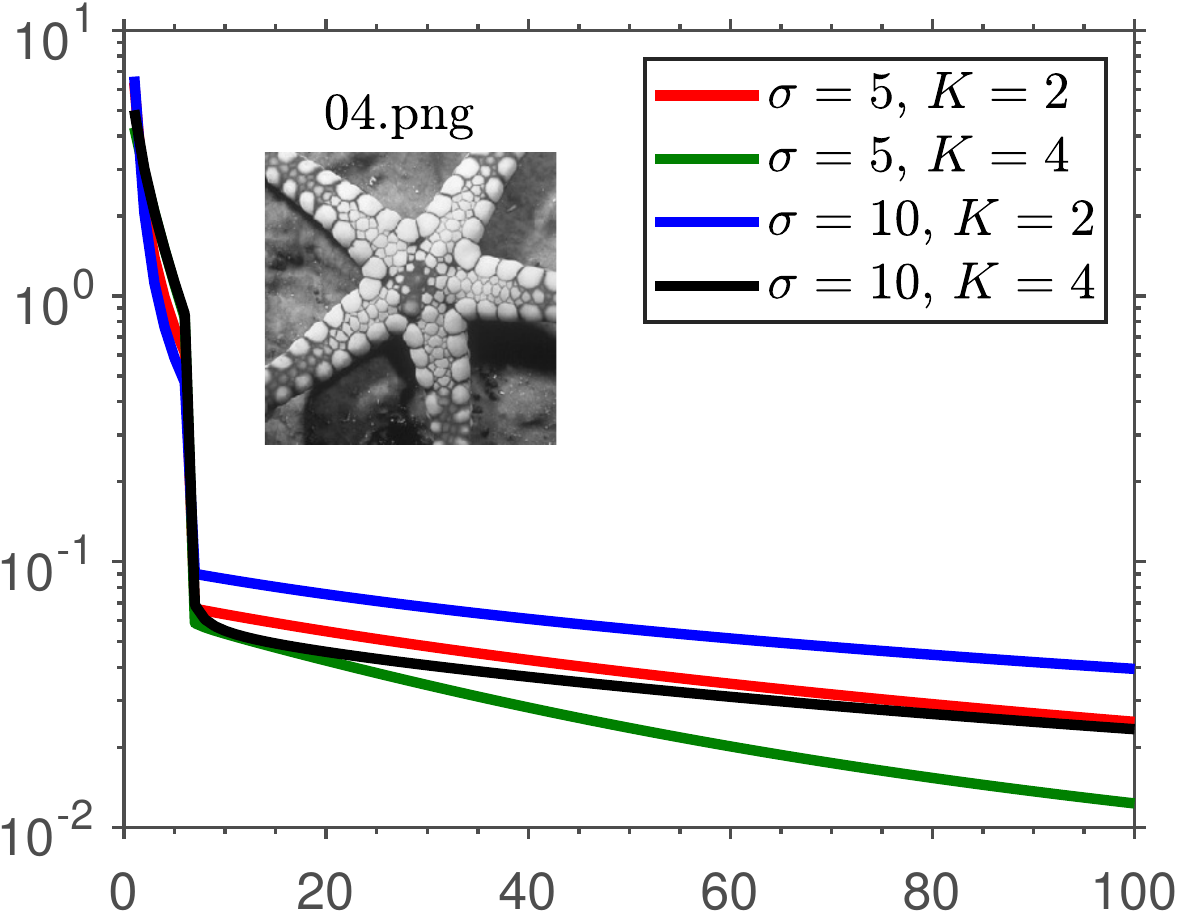}} 
\caption{Objective and residual convergence for PnP-ISTA in $(\Re^n,\inner{\cdot}{\cdot})$ for $\rho=2.5$. Top row: $f(\y_k) + \rho g_D(\y_k)$. Bottom row: $\normE{\x_{k+1} - \x_k}$ (on a logarithmic scale).} 
\label{plot:superres}
\end{figure}

In Table \ref{tab:superres}, we report  PSNR/SSIM values, averaged over the images (resized to $256 \times 256$) in the Set12 dataset \cite{zhang2017beyond}.
Note that the reconstruction quality using PnP-ISTA in $(\Re^n,\inner{\cdot}{\cdot})$ is competitive with its standard counterpart. 
In particular, we note that this behavior holds for different values of $K$ and $\sigma$.
For completeness, a visual result is shown in Fig.~\ref{fig:superres}.
We empirically verify Theorem \ref{thm:linearISTA}  for two particular images from the Set12 dataset and using different values of $K$ and $\sigma$; the results are reported in Fig.~\ref{plot:superres}.
Notice that the objective value corresponding to proposed algorithm decreases in every iteration and stabilizes.
As for the iterates $(\x_k)$, it is not possible to directly verify the convergence since the limit point $\x^\ast$ is not known.
Instead, as shown in the figure, we verify a necessary condition, namely that $\normE{\x_k - \x_{k-1}}$ decays to $0$ as $k$ increases.
%Also notice that the successive difference of iterates ${\lVert \x_k - \x_{k-1} \rVert}_2$ decay to $0$; this is necessary for iterate convergence but not sufficient though. 

\subsection{Despeckling} 

\begin{table*}
\centering
\caption{Objective and iterate convergence of PnP-ADMM in $(\Re^n,\inner{\cdot}{\cdot})$ for image despeckling ($\rho = 0.2$).}
\label{tab:despeckling_conv}
\begin{adjustbox}{max width=\linewidth}
\begin{tabular}{clccccccccccc}
\toprule
\multirow{2}{*}{Image}	& \multirow{2}{*}{Settings}	& \multicolumn{6}{c}{$f(\y_k) + \rho g_D(\y_k) - C$} 								& \multicolumn{5}{c}{$\normE{\y_{k+1} - \y_k}$}\\
\cmidrule(lr){3-8}
\cmidrule(lr){9-13}
						&							& $k=1$		& $k=10$	& $k=20$	& $k=30$	& $k=40$	& $C$		& $k=1$		& $k=10$	& $k=20$	& $k=30$	& $k=40$ \\
\midrule
\multirow{3}{*}{03.png}	& $M=5$						& $1.8 \times 10^3$	& $2.7$ & $1.6 \times 10^{-2}$ & $1.9 \times 10^{-4}$ & $3 \times 10^{-6}$ & $3.4 \times 10^6$	& $7.5 \times 10^1$ & $2.6 \times 10^{-1}$ & $1.3 \times 10^{-2}$ & $1.2 \times 10^{-3}$ & $1.4 \times 10^{-4}$ \\
						& $M=7$						& $2.7 \times 10^3$ & $8.4$ & $1.1 \times 10^{-1}$ & $2.6 \times 10^{-3}$ & $8.9 \times 10^{-5}$ & $4.7 \times 10^6$ & $6.7 \times 10^1$ & $3.2 \times 10^{-1}$ & $2.4 \times 10^{-2}$ & $3.1 \times 10^{-3}$ & $5.1 \times 10^{-4}$ \\
						& $M=10$					& $4.1 \times 10^3$ & $2.4 \times 10^1$ & $5.9 \times 10^{-1}$ & $2.6 \times 10^{-2}$ & $1.7 \times 10^{-3}$ & $6.8 \times 10^6$ & $6.3 \times 10^1$ & $3.9 \times 10^{-1}$ & $4.2 \times 10^{-2}$ & $7.3 \times 10^{-3}$ & $1.6 \times 10^{-3}$ \\
\midrule
\multirow{3}{*}{05.png}	& $M=5$						& $2.8 \times 10^3$ & $3.9$ & $1.8 \times 10^{-2}$ & $1.6 \times 10^{-4}$ & $2.0 \times 10^{-6}$ & $3.3 \times 10^6$ & $8.0 \times 10^1$ & $3.1 \times 10^{-1}$ & $1.4 \times 10^{-2}$ & $1.1 \times 10^{-3}$ & $1.1 \times 10^{-4}$ \\
						& $M=7$						& $4.1 \times 10^3$ & $1.2 \times 10^1$ & $1.2 \times 10^{-1}$ & $2.1 \times 10^{-3}$ & $5.1 \times 10^{-5}$ & $4.7 \times 10^6$ & $7.1 \times 10^1$ & $3.9 \times 10^{-1}$ & $2.7 \times 10^{-2}$ & $3.0 \times 10^{-3}$ & $4.1 \times 10^{-4}$ \\
						& $M=10$					& $6.2 \times 10^3$ & $3.5 \times 10^1$ & $6.6 \times 10^{-1}$ & $2.0 \times 10^{-2}$ & $8.2 \times 10^{-4}$ & $6.7 \times 10^6$ & $6.6 \times 10^1$ & $4.9 \times 10^{-1}$ & $4.8 \times 10^{-2}$ & $7.1 \times 10^{-3}$ & $1.3 \times 10^{-3}$ \\
\bottomrule
\end{tabular}
\end{adjustbox}
\end{table*}
In $M$-look synthetic aperture radar imaging \cite{deledalle2017mulog}, the observation model is given by 
\begin{equation}
\label{model1}
\s(\i) = \r_0(\i)\n(\i),
\end{equation}
where $\i$ is the pixel location, $\r_0 \in \mathbb{R}_{++}^{n}$ is the unknown reflectance image, and the components of $\n \in \mathbb{R}_+^{n}$ (known as \textit{speckle noise}) are i.i.d. Gamma with unit mean and variance $1/M$.
In particular, the probability density of each component  is
\begin{equation*}
p_n(t) = \frac{M^M}{\Gamma(M)} t^{M-1} \exp(-Mt) \qquad (t \geqslant 0).
\end{equation*}
The above model is based on the assumption that the measurement $\s(\i)$ is an average of $M$ independent samples of the intensity at pixel $\i$ \cite{deledalle2017mulog}. 
Letting $\tilde{\x}(\i) = \log \r_0(\i)$ and $\y(\i) = \log \s(\i)$, and taking logarithms on both sides of \eqref{model1}, we obtain the following additive model:
\begin{equation*}
\y(\i) = \tilde{\x}(\i) + \log \n(\i).
\end{equation*}
The transformed noise distribution is given by 
\begin{equation*}
p_{\log(n)}(t) = \frac{M^M}{\Gamma(M)}  \exp\left(M (t-e^{t}) \right).
\end{equation*}
The data-fidelity term corresponding to the maximum likelihood estimate of $\tilde{\x}$ is given by  \cite{bioucas2010multiplicative}:
\begin{equation*}
f(\x) = M \sum_{ \i} \Big(\x(\i) + \exp \big(\y(\i) - \x(\i)\big) \Big) + \text{constant}.
\end{equation*}
\begin{figure}[t]
\centering
\scriptsize
\stackunder{\hphantom{=}}{\stackunder{\includegraphics[width=0.24\linewidth]{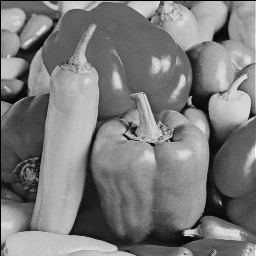}}{Ground truth.}}
\stackunder{\hphantom{=}}{\stackunder{\includegraphics[width=0.24\linewidth]{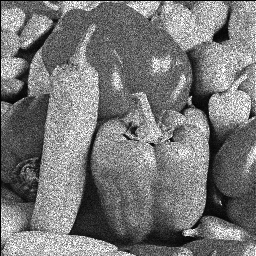}}{Observed.}}
\stackunder{PSNR = $28.77$ dB.}{\stackunder{\includegraphics[width=0.24\linewidth]{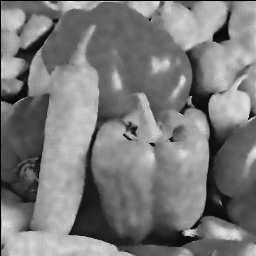}}{PnP, $(\Re^n,\innerE{\cdot}{\cdot})$.}}
\stackunder{PSNR = $28.72$ dB.}{\stackunder{\includegraphics[width=0.24\linewidth]{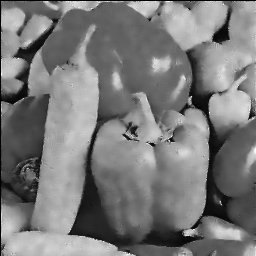}}{PnP, $(\Re^n,\inner{\cdot}{\cdot})$.}}
\caption{Image despeckling for $M=5$ using PnP-ADMM. The third and fourth images are the reconstructions from PnP-ADMM in $(\Re^n,\innerE{\cdot}{\cdot})$ and $(\Re^n,\inner{\cdot}{\cdot})$. The respective SSIM values are $0.840$ and $0.847$, while the PSNR values are reported above the images.}
\label{fig:despeckling}
\end{figure}
We use PnP-ADMM to estimate $\tilde{\x}$ using the above $f$ as the data-fidelity.
Note that the proximal operator of $f$ in the space $(\Re^n,\inner{\cdot}{\cdot})$ is given by:
\begin{equation*}
\label{eq:prox_despeckling}
\begin{aligned}
\prox_{{\rho}^{-1} f} (\u) = & \argmin_{\x \in \Re^n} \bigg\{ M \sum_{ \i} \Big[\x(\i) + \exp \big(\y(\i) - \x(\i)\big)\Big]\\
&+  \frac{\rho}{2} \left(\x - \u \right)^\top\! \D \left(\x - \u \right) \bigg\}.
\end{aligned}
\end{equation*}
The above optimization  is separable since $\D$ is diagonal. In fact, the optimization can be decoupled into  one-variable convex problems, which can be solved efficiently using Newton's method \cite{bioucas2010multiplicative}.

A visual example of despeckling of a simulated observation  is shown in Fig.~\ref{fig:despeckling} for $M=5$.
The values  of $f(\y_k) + \rho g_D(\y_k)$ and $\normE{\y_{k+1} - \y_k}$ for different values of $M$  and input images are reported in Table \ref{tab:despeckling_conv}.
It is evident that $f(\y_k) + \rho g_D(\y_k)$ converges to a stable value, whereas $\normE{\y_{k+1} - \y_k}$ decays to $0$ with iterations $k$.
These observations agree with Theorem \ref{thm:linearADMM}.
We compare the reconstruction quality with PnP-ADMM in $(\Re^n,\innerE{\cdot}{\cdot})$ by averaging the  PSNR and SSIM values over the images in the Set12 dataset (all images were resized to $256 \times 256$).
These are reported in Table \ref{tab:despeckling}.
Note that the PSNR/SSIM values are comparable for PnP-ADMM in $(\Re^n,\inner{\cdot}{\cdot})$ and $(\Re^n,\innerE{\cdot}{\cdot})$.

\begin{table}[t!]
\centering
\caption{PSNR/SSIM values for  despeckling of Set12 dataset \cite{zhang2017beyond}  using PnP-ADMM.}
\label{tab:despeckling}
\begin{tabular}{lcc}
\toprule
\diagbox{Settings}{Method} & PnP in $(\Re^n,\innerE{\cdot}{\cdot})$ & PnP in $(\Re^n,\inner{\cdot}{\cdot})$ \\
\toprule
$M=5$ & $27.65/0.809$ & $27.49/0.793$ \\
$M=10$ & $29.37/0.838$ & $29.01/0.829$ \\
\bottomrule
\end{tabular}
\end{table}

\section{Conclusion}
\label{Conclusion}

We showed that iterate and objective convergence of PnP-ISTA and PnP-ADMM can be guaranteed for a class of linear denoisers provided we work with a denoiser-specific inner product (and associated gradient and proximal operators). To the best of our knowledge, this is the first such result for PnP algorithms where simultaneous analysis of iterate and objective convergence can be carried out. Moreover, our results subsume existing convergence results for symmetric linear denoisers. Importantly, our  analysis holds for non-symmetric kernel filters like nonlocal means which is known to possess good regularization capabilities. In fact, we demonstrated this for model-based superresolution and despeckling. An interesting question arising from our work is whether the present results (and analysis) can be extended to other proximal algorithms including the accelerated variants.

\section{Appendix}
\label{sec:proofs}
In this section, we give detailed proofs of the results in Section \ref{mainresults}.
Unless specified otherwise, it should be understood that we work in $(\Re^n, \inner{\cdot}{\cdot})$, where $\inner{\cdot}{\cdot}$ is an arbitrary inner product and $\norm{\cdot}$ denotes the norm induced by $\inner{\cdot}{\cdot}$.

\subsection{Proof of Proposition~\ref{prop:averaged_l2}}
\label{subsec:proof_averaged_l2}

1) Let $\A = (1/\theta) \W - (1/\theta - 1) \I$ and $\theta \in [1/2,1)$. Since the eigenvalues of $\W$ are in $[0,1]$, the eigenvalues of  $\A$ must be in $[-1,1]$. Since $\A$ is symmetric, this means that its spectral norm $\normE{\A}$ (largest singular value)  is at most $1$. 
Hence, $\normE{\A(\x - \y)} \leqslant \normE{\A}  \normE{\x - \y} \leqslant \normE{\x - \y} $, i.e., $\A$ is non-expansive. Hence, $\W$ is $\theta$-averaged on  $(\Re^n, \innerE{\cdot}{\cdot})$.

2) Since $\W$ is row-stochastic, $\W \e = \e$, where $\e \in \Re^n$ is the all-ones vector.
Thus, $1$ is an eigenvalue of $\W$. Since $\normE{\W} = \max\{\normE{\W\x} : \normE{\x}=1\}$, we can conclude that
\begin{equation}
\label{gt1}
\normE{\W} \geqslant 1.
\end{equation}

Suppose that $D$ is averaged. Then we can show that $\W^{\top} \!\e=\e$. But this would contradict our assumption that $\W$ is not doubly stochastic, hence $D$ cannot be averaged.
Indeed, if $D$ is averaged, then it is non-expansive: $\normE{\W\x} \leqslant \normE{\x}$ for all $\x \in \Re^n$. Thus, we must have $\normE{\W} \leqslant 1$. Hence, from \eqref{gt1}, it follows that $\normE{\W} = 1$. Now, since $\normE{\W^\top} =  \normE{\W}$,
\begin{align*}
\normE{\e}^2 = {\e}^\top (\W\e) = {(\W^{\top}\e})^\top \!\e & \leqslant \normE{\W^{\top}\e} \normE{\e}\\
& \leqslant \normE{\W^\top} \normE{\e}^2 = \normE{\e}^2. 
\end{align*}
Thus, both ``$\leqslant$'' are in fact ``$=$''. In particular, for the Cauchy-Schwarz inequality, we must have $\W^{\top}\!\e=\alpha\e$, where $\alpha \in \Re$. However, $$(\alpha {\e})^\top \e = {(\W^{\top}\!\e})^\top\! \e =  {\e}^\top\! (\W \e)={\e}^\top \e.$$ 
Hence, $\alpha=1$ and  $\W^{\top} \!\e=\e$, as claimed.

\subsection{Proof of Lemma \ref{lem:key_lemma}}
To prove Lemma \ref{lem:key_lemma}, we need the following result. A proof can be found in \cite[Proposition 4.35(iii)]{bauschke2017convex}. Nevertheless, we offer a  more self-contained analysis. 
\begin{lemma}
\label{lem:averaged_bound}
Let $T : \Re^n \to \Re^n$ be $\theta$-averaged for some $\theta \in (0,1)$.
Then, for all $\x,\y \in \Re^n$,
\begin{align*}
\norm{T(\x) - T(\y)}^2 \leqslant &\norm{\x - \y}^2\\
&- \frac{1 - \theta}{\theta} \norm{(I - T)(\x) - (I - T)(\y)}^2.
\end{align*}
\end{lemma}
\begin{IEEEproof}
Let $T = (1-\theta) I + \theta S$, where $S$ is non-expansive.
Then $I - T = \theta (I - S) = \theta \Delta S$, where $\Delta S = I - S$.
We need to show that for all $\x,\y \in \Re^n$,
\begin{equation*}
\norm{T(\x) - T(\y)}^2 \leqslant \norm{\x - \y}^2 - \theta (1 - \theta) \norm{\Delta S (\x) - \Delta S (\y)}^2.
\end{equation*}
Now $T(\x) - T(\y) = (\x - \y) - \theta \big( \Delta S(\x) - \Delta S(\y) \big)$.
Hence,
\begin{align*}
& \quad \norm{T(\x) - T(\y)}^2 = \\
& \left[ \norm{\x - \y}^2 - \theta (1-\theta) \norm{\Delta S(\x) - \Delta S(\y)}^2 \right]\\
+ & \left[ \theta \norm{\Delta S(\x) - \Delta S(\y)}^2 - 2 \theta \inner{\x - \y}{\Delta S(\x) - \Delta S(\y)} \right].
\end{align*}
It suffices to show that the second term on the right is negative. Indeed, note that
\begin{align*}
&\norm{\Delta S(\x) - \Delta S(\y)}^2 - 2 \inner{\x - \y}{\Delta S(\x) - \Delta S(\y)}\\
=& \inner{\Delta S(\x) - \Delta S(\y)}{\Delta S(\x) - \Delta S(\y) - 2(\x-\y)}\\
=& \inner{\x - \y - S(\x) + S(\y)}{-\x + \y -S(\x) + S(\y)}\\
=& \norm{S(\y) - S(\x)}^2 - \norm{\y - \x}^2 \leqslant 0,
\end{align*}
since $S$ is non-expansive.
\end{IEEEproof}
%\begin{IEEEproof}
We now establish Lemma \ref{lem:key_lemma} using the above result.
Letting $S= I - T$, note that $\x_\ast \in \fix(T)$ if and only if $S(\x_\ast) = \ZE$.
Since by assumption $ \fix(T) \neq \varnothing
$, let $\x_\ast \in \fix(T)$.
Setting $\x = \x_k$ and $\y= \x_\ast$ in Lemma~\ref{lem:averaged_bound}, we have
\begin{equation}
\label{point1}
\norm{\x_{k+1} - \x_\ast}^2  \leqslant \norm{\x_{k} - \x_\ast}^2 - \frac{1 - \theta}{\theta} \norm{S(\x_k)}^2,
\end{equation}
since $S(\x_\ast) = \boldsymbol{0}$. By telescoping the sum in \eqref{point1},  we obtain
\begin{equation*}
\label{point3}
\norm{\x_{k+1} - \x_\ast}^2 \leqslant \norm{\x_{0} - \x_\ast}^2 - \frac{1 - \theta}{\theta} \sum_{j=0}^k \norm{S(\x_j)}^2.
\end{equation*}
The  quantity on the right is nonnegative for all $k \geqslant 1$. In particular, the series $\sum_{j=0}^\infty \norm{S(\x_j)}^2$ is bounded above by $ \norm{\x_{0} - \x_\ast}^2$. Thus, we must have $S(\x_k) \rightarrow \ZE$ as $k \rightarrow \infty$. 

On the other hand, it follows from \eqref{point1} that 
\begin{equation}
\label{dec}
\norm{\x_{k+1} - \x_\ast} \leqslant \norm{\x_{k} - \x_\ast} \leqslant \cdots \leqslant \norm{\x_0 - \x_\ast}. 
\end{equation}
We conclude that $(\x_k)_{k \geqslant 0}$ is bounded and thus must have a convergent subsequence, i.e., a subsequence $(\x_{k_r})_{r \geqslant 0}$ that converges to some $\hat{\x} \in \Re^n$. Since $S$ is continuous, we have
\begin{equation*}
S(\hat{\x}) = \lim_{r \rightarrow \infty} S(\x_{k_r}) = \lim_{k \rightarrow \infty} S(\x_k) = \ZE,
\end{equation*}
Hence, $\hat{\x} \in \fix(T)$. 

We are done if we can show that the original sequence $(\x_k)_{k \geqslant 0}$ converges to $\hat{\x}$.
Now, given any $\epsilon>0$, we can find $N \geqslant 1$ such that  $\norm{\x_{k_n} - \hat{\x}} < \epsilon$ for $n \geqslant N$.
Let $K = k_{N}$. Since $\hat{\x} \in \fix(T)$, it follows from \eqref{dec} that for $k \geqslant K$,
\begin{equation*}
\norm{\x_{k} -  \hat{\x}} \leqslant \norm{\x_{k_{N}} - \hat{\x}} < \epsilon.  
\end{equation*}
Since $\epsilon$ is arbitrary, we conclude that $\x_k \rightarrow \hat{\x}$ as $k \rightarrow \infty$.
%\end{IEEEproof}

\subsection{Proof of Theorem \ref{thm:PnP_ISTA_fixed}}
\label{subsec:proof_PnP-ISTA_fixed}

We need some preliminary results to establish Theorem \ref{thm:PnP_ISTA_fixed}.
\begin{lemma}
\label{lem:strictly_monotone}
Let $f : \Re^n \to \Re$ be convex and $\beta$-smooth.
Then, for any $\x,\y \in \Re^n$,
\begin{equation*}
\inner{\nabla \! f(\x) - \nabla \! f(\y)}{\x - \y} \geqslant \frac{1}{\beta} \norm{\nabla \! f(\x) - \nabla \! f(\y)}^2.
\end{equation*}
\end{lemma}

\begin{IEEEproof}
See  \cite[Theorem 2.1.5]{nesterov2013introductory}.
\end{IEEEproof}

\begin{lemma}
\label{lem:averaged_grad}
Let $f : \Re^n \to \Re$ be convex and $\beta$-smooth.
Then the operator $(I - \rho^{-1} \nabla \! f)$ is $(\beta/2 \rho)$-averaged for $\rho > \beta/2$.
\end{lemma}
\begin{IEEEproof}
Let $\theta = \beta/2 \rho \in (0,1)$, and let $S = I - (2/\beta) \nabla \! f$.
Then $I - \rho^{-1} \nabla \! f = (1-\theta) I + \theta S$.
We are done if we can show that $S$ is non-expansive.
Indeed, for any $\x,\y \in \Re^n$,
\begin{align*}
&\norm{S(\x) - S(\y)}^2 = \norm{(\x  - \y) - \frac{2}{\beta}\big( \nabla \! f(\x) - \nabla \! f(\y) \big)}^2\\
&= \norm{\x  - \y}^2 + \frac{4}{\beta} \Big[ \frac{1}{\beta} \norm{\nabla \! f(\x) - \nabla \! f(\y)}^2 \\
& \qquad \qquad \qquad \qquad - \inner{\x - \y}{\nabla \! f(\x) - \nabla \! f(\y)} \Big]\\
&\leqslant \norm{\x  - \y}^2,
\end{align*}
where the last inequality follows from Lemma~\ref{lem:strictly_monotone}.
\end{IEEEproof}

\begin{lemma}
\label{lem:averaged_comp}
If $T_1,T_2 : \Re^n \to \Re^n$ are averaged, then $T_1 \circ T_2$ is averaged.
\end{lemma}
\begin{IEEEproof}
See \cite[Proposition 4.44]{bauschke2017convex}. 
%Let $T_i = (1-\theta) I + \theta S_i, i=1,2$, where $\theta_i \in (0,1)$ and $S_i$ is non-expansive.
%Define $\theta =  \theta_1 + \theta_2 - \theta_1 \theta_2$, and
%\begin{equation*}
%S = \frac{1}{\theta} \left[ \theta_1 (1-\theta_2) S_1 + \theta_2 (1-\theta_1) S_2 + \theta_1 \theta_2 S_1 \circ S_2 \right].
%\end{equation*}
%Now $1 - \theta = (1-\theta_1)(1-\theta_2) \in (0,1)$. Moreover, from triangle inequality, we see that $S$ is non-expansive. Since $T_1 \circ T_2 = (1-\theta) I + \theta S$, it follows that $T_1 \circ T_2$ is averaged.
\end{IEEEproof}

We are now ready to prove Theorem \ref{thm:PnP_ISTA_fixed}. 
By Lemma~\ref{lem:averaged_grad}, $(I - \rho^{-1} \nabla \! f)$ is averaged if $\rho > \beta/2$.
Since $D$ is assumed to be averaged, it follows from Lemma~\ref{lem:averaged_comp} that $T_\mathrm{ISTA}$ is averaged if $\rho > \beta/2$. Finally, since $\fix(T_\mathrm{ISTA}) \neq \varnothing
$, Theorem \ref{thm:PnP_ISTA_fixed} follows immediately from Lemma~\ref{lem:key_lemma}.

\subsection{Proof of Theorem \ref{thm:PnP_ADMM_fixed}}
\label{subsec:proof_PnP-ADMM_fixed}

We first show that the PnP-ADMM iterates $(\x_k,\y_k,\z_k)_{k \geqslant 1}$ can be written in terms of a single-variable sequence $(\u_k)_{k \geqslant 1}$.
\begin{lemma}
\label{lem:T_ADMM}
Fix $\rho>0$ and $\y_0, \z_0 \in \Re^n$, and consider the sequence $(\x_{k}, \y_{k}, \z_{k})_{k \geqslant 1}$ generated by \eqref{eq:PnP_ADMM}.
Let $\u_1 = \y_1 + \z_1$, and $\u_{k} = T_{\mathrm{ADMM}}(\u_{k-1})$ for $k \geqslant 2$, where $T_{\mathrm{ADMM}}$ is given by \eqref{eq:T_ADMM}.
Then for $k \geqslant 1$, we have $\y_k = D(\u_k), \ \z_k = (I - D)(\u_k)$ and $\x_{k+1} = \prox_{\rho^{-1} f}\circ (2 D - I)(\u_k)$.
\end{lemma}
\begin{IEEEproof}
Define $\hat{\y}_k = D(\u_k), \hat{\z}_k =(I - D) (\u_k),$ and $\hat{\x}_{k+1} = \prox_{\rho^{-1} f}\circ (2 D - I)(\u_k)$.
We will inductively show that for $i \geqslant 1$,
\begin{equation}
\label{ind}
\y_i = \hat{\y}_i, \ \z_i = \hat{\z}_i  \text{ and }  \x_{i+1} = \hat{\x}_{i+1}.
\end{equation}

\underline{\textbf{Base case}}: First, note that \eqref{ind} holds for $i=1$. Indeed, from \eqref{eq:PnP-ADMM_y} and \eqref{eq:PnP-ADMM_z}, we have $\y_1 = D(\x_1 + \z_0)$ and $\z_1 = \z_0 + \x_1 - \y_1$.
Since,  by construction, $\u_1 = \y_1 + \z_1$, we have $\u_1 = \z_0 + \x_1$.
Therefore, $\y_1 = D(\x_1 + \z_0) = D(\u_1) = \hat{\y}_1$, and 
\begin{equation*}
\z_1 = \u_1 - \y_1 = \u_1 - D(\u_1) = (I - D)(\u_1) = \hat{\z}_1.
\end{equation*}
Moreover, $\y_1 - \z_1 = (2 D - I)(\u_1)$. Hence, from \eqref{eq:PnP-ADMM_x}, we get that
$\x_2 = \prox_{\rho^{-1} f} \circ (2 D - I)(\u_1) = \hat{\x}_2$.

\underline{\textbf{Induction}}: Assume that \eqref{ind} holds for $1 \leqslant i \leqslant k-1$.
We need to show that \eqref{ind} holds for $i=k$.
From the definition~of $T_{\mathrm{ADMM}}$ and the induction hypothesis, we have
\begin{align*}
\u_{k} &= T_{\mathrm{ADMM}} ( \u_{k-1} ) \\
&= \prox_{\rho^{-1} f}\circ (2 D- I)( \u_{k-1} ) + ( I - D) (\u_{k-1}) \\
&= \hat{\x}_{k} + \hat{\z}_{k-1}\\
&= \x_{k} + \z_{k-1}.
\end{align*}
Using \eqref{eq:PnP-ADMM_y}, \eqref{eq:PnP-ADMM_z} and the induction hypothesis, we see that
\begin{equation*}
\hat{\y}_{k} = D(\u_{k}) = D(\x_{k} + \z_{k-1})= \y_{k},
\end{equation*}
and
\begin{equation*}
\hat{\z}_{k} = (I - D)(\u_k)= \u_{k} - \hat{\y}_{k} = \z_{k-1} + \x_{k} - \y_{k} = \z_{k}.
\end{equation*}
Finally, note that $(2D -I)(\u_k) = \y_k - \z_k$.
Thus, from \eqref{eq:PnP-ADMM_x} and the definition~of $\hat{\x}_{k+1}$, we have
\begin{equation*}
\hat{\x}_{k+1} = \prox_{\rho^{-1} f} (\y_k - \z_k) = \x_{k+1}.
\end{equation*}
This completes the induction and the proof of Lemma \ref{lem:T_ADMM}.
\end{IEEEproof}

Next, we need a standard result about proximal operators.
\begin{lemma}
\label{lem:averaged_prox}
Let $f : \Re^n \to \overline{\Re}$ be closed, proper and convex.
Then $\prox_f$ is $(1/2)$-averaged.
\end{lemma}
\begin{IEEEproof}
See \cite[Proposition 2.28, Proposition 4.4]{bauschke2017convex}.
\end{IEEEproof}
Finally, we need the following property of averaged operators.
\begin{lemma}
\label{lem:averaged_property}
If $T: \Re^n \to \Re^n$ is $\theta$-averaged, then $T$ is $\bar{\theta}$-averaged for all $\bar{\theta} \in [\theta, 1)$. 
\end{lemma}
\begin{IEEEproof}
Let $T = (1-\theta) I + \theta S$, where $\theta \in (0,1)$ and $S$ is non-expansive. Fix $\bar{\theta} \in [\theta, 1)$ and define $\gamma =  \theta/ \bar{\theta}$. Then $T = (1-\bar{\theta}) I + \bar{\theta} \bar{S}$, where $\bar{S} = (1-\gamma) I + \gamma S$. Since $\bar{S}$ is non-expansive by triangle inequality and $\bar{\theta} \in [\theta, 1)$ was arbitrary, it follows that $T$ is averaged for all $\bar{\theta} \in [\theta, 1)$. 
\end{IEEEproof}

Using the above results, we can now establish Theorem \ref{thm:PnP_ADMM_fixed}. Note that:
\begin{itemize}
\item $2 \prox_{\rho^{-1} f} - I$ is non-expansive since $\prox_{\rho^{-1} f}$ is $(1/2)$-averaged (Lemma~\ref{lem:averaged_prox}).
\item $2 D - I$ is non-expansive since $D$ is $(1/2)$-averaged by Lemma \ref{lem:averaged_property}.
\item Thus, $(2 \prox_{\rho^{-1} f} - I) \circ (2D - I)$ is non-expansive. Hence, $T_{\mathrm{ADMM}}$ is $(1/2)$-averaged.  
\item Since $\fix(T_{\mathrm{ADMM}}) \neq \varnothing$, Lemma~\ref{lem:key_lemma} implies that the sequence ${(\u_{k})}_{k \geqslant 1}$ in Proposition~\ref{lem:T_ADMM} converges to some $\u^{\ast} \in \fix(T_{\mathrm{ADMM}})$. 
\item Since $D$ and $\prox_{\rho^{-1} f}$ are continuous, we have
\begin{align*}
\x_k  &\longrightarrow \prox_{{\rho}^{-1} f}\circ \left(2 D - I \right)(\u^\ast),\\
\y_k &\longrightarrow  D(\u^\ast)  \ \text{ and } \ \z_k  \longrightarrow (I - D)(\u^\ast). 
\end{align*}
\end{itemize}
This completes the proof of Theorem \ref{thm:PnP_ADMM_fixed}.

\subsection{Proof of Theorem \ref{thm:averaged_combined}}
\label{subsec:proof_averaged_combined}
(a) We use the following result from \cite[Corollary 3.4]{combettes2018monotone}:
A linear operator $D$ on $(\Re^n,\inner{\cdot}{\cdot})$ is the proximal operator of a closed proper convex function if the following conditions are met: (i) $D$ is self-adjoint; (ii) $\inner{\x}{D(\x)} \geqslant 0$ for all $\x \in \Re^n$; and (iii) $\norm{D} \leqslant 1$, where $\norm{D}$ is the operator norm of $D$.

We verify that the above three conditions are satisfied by any $D \in \Lin$.
Let $\W = \V \A \V^{-1}$ be an eigendecomposition of $\W$, where the diagonal matrix $\A$ contains the eigenvalues of $\W$
For (i), we need to show that $\inner{\x}{D(\y)} = \inner{D(\x)}{\y}$ for all $\x,\y \in \Re^n$.
This is straightforward from definition \eqref{inner} by writing $\W = \V \A \V^{-1}$.
Similarly, (ii) follows using the fact that the entries of $\A$ are nonnegative.
To prove (iii), note that by definition $\norm{D} = \max \{ \norm{\W \x} : \norm{\x} = 1 \}$.
Let $\y = \V^{-1} \x$. Then $\norm{\x} = \norm{\y}_2$. Moreover,  since $\W = \V \A \V^{-1}$, we have $\norm{\W \x} =  \normE{\A \y}$. Hence, 
\begin{equation*}
\norm{D} = \max \left\{ \normE{\A \y} : \normE{\y} = 1 \right \}  \leqslant 1,
\end{equation*}
since the diagonal entries of $\A$ are in $[0,1]$.

(b) From part (a), we have
\begin{equation}
\label{den_eqn}
\W \x = \prox_{g_{D}}(\x), \qquad (\x \in \Re^n)
\end{equation}
where $g_D : \Re^n \to \Rc$ is   closed proper and convex. Hence (b) follows from Lemma \ref{lem:averaged_prox} and Lemma \ref{lem:averaged_property}. 

(c) From \eqref{den_eqn}, $\prox_{g_{D}}$ maps $\Re^n$ onto $\mcR(\W)$.
Moreover, since $g_D$ is proper and from \eqref{den_eqn}, 
%for any $\x \in \Re^n$, $g_D \big( \prox_{g_{D}}(\x) \big) < \infty$ since $g_D$ is proper.
%Hence, 
$g_D$ is real-valued on $\mcR(\W)$.
%Hence, $\mcR(\W) \subseteq \dom (g_{D})$, i.e., $g_{D}$ is real-valued on $\mcR(\W)$.
To prove that $g_{D}$ restricted to $\mcR(\W)$ is continuous, note that $\mcR(\W)$ is a subspace of $\Re^n$.
If $\dim \mcR(\W) = 0$, then continuity follows trivially.
Hence, assume that $\dim \mcR(\W) = p \geqslant 1$.
Let $\varphi: \Re^p \to \mcR(\W)$ be a linear isomorphism, and define $h = g_{D} \circ \varphi$.
Since $g_D$ is convex and $\varphi$ is linear, $h$ is convex.
Since real-valued convex functions are continuous (see \cite[Corollary 8.40]{bauschke2017convex}), $h$ is continuous.
Further, since $\varphi^{-1}$ is continuous, it follows that $g_{D}|_{\mcR(\W)} = h \circ \varphi^{-1}$ is continuous.

(d) Since $D = \prox_{g_{D}}$, we can replace the denoiser $D$ by $ \prox_{g_{D}}$ in \eqref{eq:PnP_ISTA} and \eqref{eq:PnP-ADMM_y}. Then, the updates in PnP-ISTA (resp. PnP-ADMM) are exactly that of ISTA (resp. ADMM) applied to optimization problem \eqref{eq:mainoptlinear}.

%\end{IEEEproof}

\subsection{Proof of Theorem \ref{thm:linearISTA}}
\label{subsec:proof_pnp_linear_conv}
To establish this theorem, we will require the concept of subdifferential of a convex function. 
\begin{definition}
\label{def:subdifferential}
Let $f : \Re^n \to \Rc$ be proper and convex.
The subdifferential of $f$ at $\x \in \Re^n$, denoted by $\partial f(\x)$, is defined to be the set of all vectors $\v \in \Re^n$ such that
\begin{equation*}
f(\y) \geqslant f(\x) + \inner{\v}{\y - \x} \qquad (\y \in \Re^n).
\end{equation*}
\end{definition}
As with the rest of the discussion, we assume $\inner{\cdot}{\cdot}$ to be an arbitrary (but fixed) inner product on $\Re^n$.

We need the following properties of subdifferential.
\begin{lemma}
\label{lem:subgradient_sum}
Let $f : \Re^n \to \Re$ be convex and $g : \Re^n \to \Rc$ be closed proper and convex. 
%Let $\inner{\cdot}{\cdot}$ be an arbitrary inner product on $\Re^n$. 
Then
\begin{enumerate}[(a)]
\item $\x^\ast$ minimizes $h(\x) = g(\x) + (\rho/2) \norm{\x - \bar{\x}}^2$, where $\rho > 0$ and $\bar{\x} \in \Re^n$, if and only if $\ZE \in \partial g(\x^\ast) + \rho (\x^\ast - \bar{\x})$.
\item $\x^\ast$ minimizes  $f + g$ if and only if $\ZE \in \partial f(\x^\ast) + \partial g(\x^\ast)$.
If $f$ is differentiable, then $\ZE \in \nabla \! f(\x^\ast) + \partial g(\x^\ast)$.
\end{enumerate}
\end{lemma}
\begin{IEEEproof}
For (a), see \cite[Th. 16.3 and Example 16.43]{bauschke2017convex}.
For (b), see \cite[Th. 16.3 and Corollary 16.48]{bauschke2017convex}.
\end{IEEEproof}
\begin{proposition}
Let $f: \Re^n \to \Re$ be convex and differentiable, and let $g: \Re^n \to \overline{\Re}$ be closed proper and convex.
Then $\x^\ast \in \mathrm{fix}({T}_{\mathrm{ISTA}})$ where $D = \prox_{\rho^{-1} g}$ if and only if $\x^\ast$ is a  minimizer of $f + g$.  
\label{prop:FBSobjequivalence} 
\end{proposition}
\begin{IEEEproof}
Let $\x^\ast \in \mathrm{fix}(T_{\mathrm{ISTA}})$, where $T_{\mathrm{ISTA}} = \prox_{{\rho^{-1}} g} \circ (I - \rho^{-1} \nabla \! f)$. Let $\y^\ast = (I - \rho^{-1} \nabla \! f) (\x^\ast)$.
Then $ \nabla f (\x^\ast) = \rho (\x^\ast - \y^\ast)$. Moreover, 
\begin{equation*}
\x^\ast = T_{\mathrm{ISTA}}(\x^\ast) = \prox_{{\rho^{-1}} g}(\y^\ast).
\end{equation*}
That is,
\begin{equation*}
\x^\ast =  \argmin_{\y \in \Re^n} \ g(\y)+ \frac{\rho}{2} \norm{\y - \y^\ast}^2.
\end{equation*}
Therefore, by Lemma~\ref{lem:subgradient_sum}, $\ZE \in \partial g (\x^\ast) + \rho (\x^\ast - \y^\ast)$.
Thus,
\begin{equation*}
\ZE \in  \partial g (\x^\ast) + \rho (\x^\ast - \y^\ast) = \partial g (\x^\ast) + \nabla \! f (\x^\ast). 
\end{equation*} 
Using Lemma~\ref{lem:subgradient_sum}, we conclude that $\x^\ast$ is a minimizer of $f + g$. 
The other direction can be established by reversing the above steps.
\end{IEEEproof}

Using the above results, we can now arrive at Theorem \ref{thm:linearISTA}. Note that since $D \in \Lin$, $D$ is averaged by Proposition~\ref{thm:averaged_combined}. Moreover, since \eqref{eq:mainoptlinear} has a minimizer, $\fix(T_\mathrm{ISTA}) \neq \varnothing$  by Proposition~\ref{prop:FBSobjequivalence}. Therefore, from Theorem \ref{thm:PnP_ISTA_fixed} and Proposition~\ref{prop:FBSobjequivalence}, we can conclude that ${(\x_k)}_{k \geqslant 0} \to \x^\ast$, where $\x^\ast$ is a minimizer of  \eqref{eq:mainoptlinear}. In particular, since $D = \prox_{g_D}$,
\begin{equation*}
p^\ast=f(\x^\ast) + \rho g_{D}(\x^\ast) .
\end{equation*}
Now, $\x^\ast \in \fix(T_\mathrm{ISTA}) \subseteq \mcR(\W)$ and $\x_k \in \mcR(\W)$ for all $k$.
On the other hand, $f$ is continuous on $\Re^n$ and $g_{D}$ is continuous on $\mcR(\W)$ by Theorem \ref{thm:averaged_combined}. Therefore, by continuity, we have
\begin{equation*}
f(\x_k) + \rho g_{D}(\x_k) \longrightarrow f(\x^\ast) + \rho g_{D}(\x^\ast) = p^\ast.
\end{equation*}
This completes the proof of Theorem  \ref{thm:linearISTA}.

\subsection{Proof of Theorem \ref{thm:linearADMM}}
\label{subsec:proof_pnp_admm_linear}
We first prove a result similar to Proposition \ref{prop:FBSobjequivalence} that relates $\mathrm{fix}(T_{\mathrm{ADMM}})$ to the minimizers of \eqref{eq:mainoptlinear}.
\begin{proposition}
Let $f: \Re^n \to \Re$ be convex and $g: \Re^n \to \overline{\Re}$ be closed, proper and convex.
Then $\x^\ast \in \mathrm{fix}({T}_{\mathrm{ADMM}})$ with $D = \prox_{\rho^{-1} g}$ if and only if $D(\x^\ast)$ is a minimizer of $f + g$.. 
\label{prop:ADMMobjequivalence} 
\end{proposition}
\begin{IEEEproof}
Suppose $\x^\ast \in \mathrm{fix}(T_{\mathrm{ADMM}})$, i.e.,
\begin{equation*}
\x^\ast = \left(2 \prox_{\rho^{-1} f} - I\right) \circ \left(2\prox_{\rho^{-1} g} - I\right)(\x^\ast).
\end{equation*}
Let $\y^\ast = \left(2\prox_{\rho^{-1} g} - I\right) (\x^\ast)$. Note that
\begin{equation}
\label{valassign}
\prox_{\rho^{-1} g} (\x^\ast) = \frac{\x^\ast + \y^\ast}{2} = \prox_{\rho^{-1} f} (\y^\ast). 
\end{equation} 
From Lemma~\ref{lem:subgradient_sum}, we have
\begin{align*}
\ZE &\in \partial g \Big(\frac{\x^\ast + \y^\ast}{2}\Big) + \frac{\rho}{2} (\y^\ast - \x^\ast), \ \text{and}\\
\ZE &\in \partial f \Big(\frac{\x^\ast + \y^\ast}{2}\Big) + \frac{\rho}{2} (\x^\ast - \y^\ast).
\end{align*}
Adding these inclusions and using \eqref{valassign}, we obtain
\begin{equation*}
\ZE \in \partial f \left(\prox_{\rho^{-1} g}(\x^\ast) \right) + \partial g \left(\prox_{\rho^{-1} g}(\x^\ast)\right). 
\end{equation*}
Since $D=\prox_{\rho^{-1} g}$, we conclude from Lemma~\ref{lem:subgradient_sum} that $D(\x^\ast)$ is a minimizer of $f + g$. 

Conversely, let $\y^\ast = \prox_{\rho^{-1} g}(\x^\ast)$.
From Lemma~\ref{lem:subgradient_sum}, we have $\ZE \in \partial g (\y^\ast) + \rho (\y^\ast - \x^\ast)$.
Moreover, since $\y^\ast = D(\x^\ast)$, $\y^\ast$ is a minimizer of $f + g$.
Using Lemma~\ref{lem:subgradient_sum} once more, we have $\ZE \in \partial f (\y^\ast) + \partial g (\y^\ast)$.
Combining these, we obtain
\begin{equation*}
\ZE \in \partial f (\y^\ast) + \rho (\x^\ast - \y^\ast),
\end{equation*}
i.e., $\y^\ast = \prox_{ \rho^{-1} f}(2 \y^\ast - \x^\ast)$, or equivalently,
\begin{equation}
\label{prop7_eqn3}
\x^\ast = \left(2 \prox_{\rho^{-1} f} - I \right) (2 \y^\ast - \x^\ast).
\end{equation}
However, $2 \y^\ast - \x^\ast = \left(2 \prox_{\rho^{-1} g} - I \right) (\x^\ast)$. Therefore, from \eqref{prop7_eqn3}, we have
\begin{equation*}
\x^\ast = \left(2 \prox_{\rho^{-1} f} - I\right) \circ \left(2\prox_{\rho^{-1} g} - I\right)(\x^\ast),
\end{equation*}
i.e., $\x^\ast \in \fix({T}_{\mathrm{ADMM}})$ with $D = \prox_{\rho^{-1}g}$.
\end{IEEEproof}

We can now establish Theorem \ref{thm:linearADMM}.
Note that since $D \in \Lin$, $D$ is $(1/2)$-averaged by Theorem~\ref{thm:averaged_combined}. Also, since \eqref{eq:mainoptlinear} has a minimizer, $\fix(T_\mathrm{ADMM}) \neq \varnothing$ by Proposition~\ref{prop:ADMMobjequivalence}. Using Theorem \ref{thm:PnP_ADMM_fixed}, we can conclude that there exists $\u^\ast \in \fix(T_\mathrm{ADMM})$ such that ${(\y_k)}_{k \geqslant 0}$ converges to $\y^\ast = \prox_{g_{D}}(\u^\ast)$.
Hence by Proposition~\ref{prop:ADMMobjequivalence}, $\y^\ast$ is a minimizer of \eqref{eq:mainoptlinear}, i.e., $f(\y^\ast) + \rho g_D(\y^\ast) = p^\ast$.
Finally, note that $\y^\ast = D(\u^\ast) \in \mcR(\W)$ and
$(\y_k) \subseteq \mcR(\W)$ for all $k$.
 Since $g_D$ is continuous on $\mcR(\W)$ by Theorem \ref{thm:averaged_combined} and $f$ is continuous, $f(\y_k) + \rho g_{D}(\y_k) \longrightarrow f(\y^\ast) + \rho g_{D}(\y^\ast) = p^\ast$. This completes the proof of Theorem \ref{thm:linearADMM}.

\subsection{Proof of Proposition~\ref{prop:neighborfilter}}
\label{subsec:proof_neighborfilter}
%\begin{IEEEproof}
Since the kernel $\phi$ is symmetric and positive definite, it follows from Definition \ref{def:pdkernel} that the kernel matrix $\K$ in \eqref{kernelmatrix} is symmetric and positive semidefinite.
Thus, the matrix $\W = \D^{-1} \K$ is nonnegative and stochastic.
Therefore, $\W$ has a complete set of eigenvectors with eigenvalues in $[0,1]$ \cite{Singer2009_diffusion_nonlocal}, which implies that $D \in \Lin$.

For the second part, we apply Theorem \ref{thm:averaged_combined}(b). This tells us that $D$ is $\theta$-averaged on $(\Re^n,\inner{\cdot}{\cdot})$ for every $\theta \in [1/2,1)$, where $\inner{\cdot}{\cdot}$ is given by \eqref{inner}. In particular, let  $\B = \D^{-\frac{1}{2}} \K \D^{-\frac{1}{2}}$ and and $\B = \U \Sigma\U^\top$ be its eigendecomposition. We can then write $\W = \D^{-\frac{1}{2}} \B \D^{\frac{1}{2}} = \V \Sigma \V^{-1}$, where $\V = \D^{-\frac{1}{2}}\U$ is an eigen matrix of $\W$.
Using this eigen matrix, $\inner{\cdot}{\cdot}$ is given by
\begin{equation*}
\inner{\x}{\y} = \innerE{\V^{-1}\x}{\V^{-1}\y} = \x^\top\! \D \y,
\end{equation*}
where we have used the fact that $\U^\top \U = \I$. This completes the proof of Proposition~\ref{prop:neighborfilter}. 

\subsection{Proof of Theorem \ref{thm:NLMkernelfilt}}

We will first establish that $\K$ is positive definite using the following result \cite[Chapter 13, Lemma 6]{cheney2009course}.  In this part, we use $\iota$ to denote $\sqrt{-1}$.

\begin{proposition}
\label{propcomplex}
Let $\x_1, \ldots, \x_n$ be distinct points in $\Re^n$ and let $c_1, \ldots, c_n \in \Re$ be such that not all are zero. Then the function $h: \Re^n \to \mathbb{C}$ defined as 
\begin{equation*}
h(\y) = \sum_{i=1}^n c_i \exp \left( \iota {\inner{\y}{\x_k}}_2 \right)
\end{equation*}
is non-zero almost everywhere. 
\end{proposition}

Since $\K$ is symmetric, we need to show that $\c^\top \K \c > 0$ for all non-zero $\c \in \Re^N$.
Let $c_{\s}$ denote the element of $\c$ corresponding to the spatial location $\s \in \Omega \subset \mathbb{Z}^2$, where $\Omega$ is the support of the image.
From \eqref{nlmkernel}, we have
\begin{align*}
\c^\top \K \c &= \sum_{\s \in \Omega} \sum_{\t \in \Omega} c_{\s} c_{\t} \, \Lambda(\s - \t) \kappa(\P_{\s} - \P_{\t}) \\ 
&= \sum_{\s \in \Omega} \sum_{\t \in \Omega} \Bigg\{ c_{\s} c_{\t} \, \int_{\Re^2} \hat{\Lambda}(\boldsymbol{\omega}_1) e^{\iota \boldsymbol{\omega}_1^\top (\s - \t)} d \boldsymbol{\omega}_1\\
& \hphantom{=====} \times \int_{\Re^{P}} \hat{\kappa}(\boldsymbol{\omega}_2) e^{\iota \boldsymbol{\omega}_2^\top (\P_{\s} - \P_{\t})} d \boldsymbol{\omega}_2 \Bigg\}
%&= \sum_{p=1}^N \sum_{q=1}^N c_p c_q \int_{\Re^M} \hat{\phi}(\boldsymbol{\omega}) \exp \left(\iota \boldsymbol{\omega}^\top\! (\z_p - \z_q) \right) d\boldsymbol{\omega}. 
\end{align*}
where $\hat{\Lambda}$ and  $\hat{\kappa}$ are the Fourier transforms of $\Lambda$ and $\kappa$, and $P$ is the patch size.
%In particular, we can verify that $\hat{\phi}$ is positive and continuous.
Let $\boldsymbol{\omega} = (\boldsymbol{\omega}_1,\boldsymbol{\omega}_2) \in \Re^{2+P}$ and $\boldsymbol{\zeta}_{\s} = (\s,\P_{\s}) \in \Re^{2+P}$.
Switching the sums and the integrals, and using Fubini's theorem, we get
\begin{equation}
\label{integ}
\c^\top \K \c = \int_{\Re^{2+P}} \hat{\Lambda}(\boldsymbol{\omega}_1) \hat{\kappa}(\boldsymbol{\omega}_2) |h(\boldsymbol{\omega})|^2 d\boldsymbol{\omega},
\end{equation}
where $h(\boldsymbol{\omega}) = \sum_{\s \in \Omega} c_{\s} \exp(\iota \boldsymbol{\omega}^\top\! \boldsymbol{\zeta}_{\s})$.
%\begin{equation}
%\label{integ}
%\c^\top \K \c = \int_{\Re^M} \hat{\phi}(\boldsymbol{\omega})  |h(\boldsymbol{\omega})|^2 d\boldsymbol{\omega},
%\end{equation}
%where $h(\boldsymbol{\omega}) = \sum_{i=1}^N c_i \exp (\iota \boldsymbol{\omega}^\top\! \z_i)$. 
%Note that $\hat{\phi}(\boldsymbol{\omega})>0$ almost everywhere as fourier transform of $\kappa(\cdot)$ is positive and fourier transform of $\Lambda(\cdot)$ is separable $\mathrm{sinc}^2(\cdot)$ function which is positive almost everywhere. 
Now, note that the points $\boldsymbol{\zeta}_{\s}$ are distinct for all $\s \in \Omega$.
% since the pixel locations $\s \in \Omega$ are different for different pixels. 
Therefore, we can conclude from Proposition \ref{propcomplex} that $h(\boldsymbol{\omega})$ is non-zero almost everywhere on $\Re^{2+P}$.
Moreover, since $\hat{\Lambda}$ is positive almost everywhere on $\Re^2$ and $\hat{\kappa}$ is positive everywhere on $\Re^P$, the function $\boldsymbol{\omega} \mapsto \hat{\Lambda}(\boldsymbol{\omega}_1) \hat{\kappa}(\boldsymbol{\omega}_2)$ is positive almost everywhere on $\Re^{2+P}$.
Hence, the integral \eqref{integ} is positive.
This establishes the claim that $\K$ is positive definite.

Next, note that since $\K$ is invertible, we can define
\begin{equation}
\label{gD}
g_{D}(\x) = \frac{1}{2} \x^\top \D (\K^{-1} \D - \I) \x.
\end{equation}
Now, $\D (\K^{-1} \D - \I) = \D \K^{-1} \D - \D$ is symmetric and we can verify that its eigenvalues are nonnegative.
Hence, \eqref{gD} is convex.
Furthermore, note that $\W = \D^{-1} \K$ and $g_{D}(\x) = (1/2) \inner{\x}{(\K^{-1} \D - \I)\x}$, where $\inner{\x}{\x} = \x^\top\! \D\x$ is the inner product in Proposition \ref{thm:NLMkernelfilt}. We claim that $\W\x = \prox_{g_{D}} (\x)$, i.e.,
\begin{equation}
\label{obj}
\W\x  = \argmin_{\y \in \Re^n}\  \frac{1}{2} \inner{\y}{(\K^{-1} \D - \I)\y} + \frac{1}{2} \norm{\y - \x}^2,
\end{equation}
where $\norm{\cdot}$ is induced by $\inner{\cdot}{\cdot}$. This follows from the observation that the derivative of the objective in \eqref{obj}  is zero at $\y=\W\x$.

\section*{Acknowledgements}
We thank the Associate Editor and the anonymous reviewers for examining the manuscript in detail and for their comments and suggestions.

% -----------------------------------------------------
\bibliographystyle{IEEEtran}
\bibliography{refs}
% -----------------------------------------------------

\end{document}